\documentclass[final,3p,times]{elsarticle}

\usepackage{amsmath,amssymb,amscd}
\usepackage{enumitem}
\usepackage{lineno}
\modulolinenumbers[5]
\usepackage{graphicx}
\usepackage{epic,eepic}
\usepackage{xcolor}
\usepackage{subfigure}
\usepackage{algorithm}
\usepackage{algorithmic}
\usepackage{diagbox}
\usepackage{url}
\usepackage{hyperref}
\usepackage{tikz}
\usepackage{xparse}
\usepackage{booktabs,cellspace}
\usepackage{siunitx}
\usepackage{caption} 
\usepackage{pgfplots}
\usepackage{pgfplotstable}
\usepgfplotslibrary{patchplots}

\biboptions{sort&compress}
\usetikzlibrary{matrix,backgrounds}
\pgfdeclarelayer{myback}
\pgfsetlayers{myback,background,main}

\tikzset{mycolor/.style = {line width=1bp,color=#1}}%
\tikzset{myfillcolor/.style = {draw,fill=#1}}%

\NewDocumentCommand{\highlight}{O{blue!40} m m}{%
\draw[mycolor=#1] (#2.north west)rectangle (#3.south east);
}

\NewDocumentCommand{\fhighlight}{O{blue!40} m m}{%
\draw[myfillcolor=#1] (#2.north west)rectangle (#3.south east);
}

\newcolumntype{j}{>{\hsize=.25\hsize\centering\arraybackslash}c}

\newcommand{\set}[1]{\mathcal{#1}}

\def\bu{{\boldsymbol u}}

\newcommand{\bx}{{\boldsymbol x}}

\DeclareMathOperator*{\argmax}{arg\,max}
\newdefinition{rmk}{Remark}
\newdefinition{lemma}{Lemma}
\newdefinition{definition}{Definition}
\newdefinition{proposition}{Proposition}
\newproof{proof}{Proof}

\definecolor{butter1}{rgb}{0.988,0.914,0.310}
\definecolor{chocolate1}{rgb}{0.914,0.725,0.431}
\definecolor{chameleon1}{rgb}{0.541,0.886,0.204}
\definecolor{skyblue1}{rgb}{0.247,0.524,0.912}
\definecolor{applegreen}{rgb}{0.55, 0.71, 0.0}
\definecolor{blue-green}{rgb}{0.0, 0.87, 0.87}
\definecolor{plum1}{rgb}{0.678,0.498,0.659}
\definecolor{scarletred1}{rgb}{0.937,0.161,0.161}

\journal{My journal}

\begin{document}
\title{Comparison of nonlinear field-split preconditioners for two-phase flow in heterogeneous porous media}

\begin{frontmatter}
\author[Stanford]{Mamadou N'diaye\corref{cor1}}
\ead{mndiaye@uphf.fr/ndiaye.mamadou24@gmail.com}
\author[Total]{Fran\c cois P. Hamon}
\ead{francois.hamon@totalenergies.com}
\author[Stanford]{Hamdi A. Tchelepi}
\ead{tchelepi@stanford.edu}

\cortext[cor1]{Corresponding author.}
\address[Stanford]{Energy Resources Engineering, Stanford University, Stanford, United States}
\address[Total]{TotalEnergies E\&P Research \& Technology USA, LLC, United States}

\begin{abstract}
This work focuses on the development of a two-step field-split nonlinear preconditioner to accelerate the convergence of two-phase flow and transport in heterogeneous porous media.
We propose a field-split algorithm named Field-Split Multiplicative Schwarz Newton (FSMSN), consisting in two steps: first, we apply a preconditioning step to update pressure and saturations nonlinearly by solving approximately two subproblems in a sequential fashion; then, we apply a global step relying on a Newton update obtained by linearizing the system at the preconditioned state.
Using challenging test cases, FSMSN is compared to an existing field-split preconditioner,  Multiplicative Schwarz Preconditioned for Inexact Newton (MSPIN), and to standard solution strategies such as the Sequential Fully Implicit (SFI) method or the Fully Implicit Method (FIM). 
The comparison highlights the impact of the upwinding scheme in the algorithmic performance of the preconditioners and the importance of the dynamic adaptation of the subproblem tolerance in the preconditioning step. 
Our results demonstrate that the two-step nonlinear preconditioning approach---and in particular, FSMSN---results in a faster outer-loop convergence than with the SFI and FIM methods.
The impact of the preconditioners on computational performance--i.e., measured by wall-clock time--will be studied in a subsequent publication.
\end{abstract}


\end{frontmatter}

%%%%%%%%%%%%%%%%%%%%%%%%%%%%%%%%%%%%%%%%%%%%%%%%%%%%%%%%%%%%%%%%%%%%%%%%%%%%%%%%%%%%%%% 
%                           Introduction
%%%%%%%%%%%%%%%%%%%%%%%%%%%%%%%%%%%%%%%%%%%%%%%%%%%%%%%%%%%%%%%%%%%%%%%%%%%%%%%%%%%%%%%

\section{Introduction}\label{sec:intro}

The numerical simulation of multiphase flow and transport in geological porous media requires solving complex partial differential equations (PDEs) with highly nonlinear saturation-dependent coefficients \cite{aziz1979petroleum,peacemanfundamentals}. 
In most subsurface applications, the heterogeneity of the porous medium generates large spatial variations in the flow regimes with high velocities near the wells and in high-permeability regions. 
This leads to severe stability constraints on the time step size in explicit discretization schemes.
Therefore, unconditionally stable implicit schemes are often the temporal discretization method of choice for porous media flow problems.

However, solving the nonlinear systems resulting from an implicit discretization of the PDEs is challenging and often represents most of the computational cost of the simulations.
In the Fully Implicit Method (FIM), all the degrees of freedom of the system--typically, pressure, saturations, and compositions--are updated simultaneously using Newton's method with damping \cite{deuflhard2011newton}.
This requires solving large, ill-conditioned linear systems at each nonlinear iteration which is computationally expensive. 
In addition, for highly nonlinear problems and/or large time steps, Newton's method often fails to converge, in which case the time step is restarted from the previous converged state with a reduced size. 
To avoid these convergence failures, globalization methods have been developed to enlarge the convergence radius of Newton's method. 
Popular damping methods for the Newton updates include saturation chopping based on heuristics \cite{younis2011modern} or relying on the structure of the fractional flow function \cite{jenny2009,wang2013,li2015,NonlinearSolverForThreePhaseTransportProblemsBasedOnApproximateTrustRegions_Moyner}.
In the Sequential Fully Implicit (SFI) method \cite{aMassConservativSequentialImplicitMultiscaleMethodForIsothermalEquationOfStateCompositionalProblems_MoynerTchelepi,aMultiscaleRestrictionSmoothedBasisMethodForCompositionalModels_MoynerTchelepi,ConsistentUpwindingForSFImultiscaleCompositionalSimulation_MoncorgeMoynerTchelepi,SuccessfulApplicationOfMultiscaleMethodsInArealReservoirSimulatorEnvironment_LieMoynerNatvig_et_al,SFIformulationForCompositionalSimulationUsingNaturalVariables_MoncorgeTchelepiJenny}, the system is decomposed into two subproblems--namely, flow and transport--to avoid solving large coupled linear systems involving all the degrees of freedom. 
In this approach, the subproblems are solved sequentially using specialized nonlinear solvers until convergence of the outer loop is attained. 
However, this decoupled approach often suffers from slow outer loop convergence for tightly coupled physics (e.g., in the presence of strong buoyancy and capillary effects) and may require the use of dedicated convergence accelerators for challenging problems \cite{NonlinearAccelerationOfSFI_Jiamin,CouplingStrengthCriteriaForSequentialImplicitFormulations_FrancMoynerTchelepi,li2021sequential,aStableMultiphaseNonlinearTransportSolverWithHybridUpwindDiscretizationInMultiscaleReservoirSimulator_WatanabeLiBratvedtLeeNatvig}.

In recent years, advanced nonlinear strategies have been proposed to overcome the limitations of Newton-based FIM and of SFI for multiphase flow in porous media.
Although a comprehensive review is out of the scope of the present work, we mention some of the recently presented nonlinear solution strategies.
They include homotopy methods \cite{deuflhard2011newton,younis2010adaptively,jiang2018dissipation}, in which a sequence of prediction-correction steps is used to follow a suitably parameterized homotopy path leading to the solution.
Homotopy methods are robust for large time steps and can achieve, in some cases, unconditional nonlinear convergence.
Using a different approach, ordering-based methods \cite{natvig2008fast,PotentialBasedReducedNewtonAlgorithmForNonlinearMultiphaseFlowInPorousMedia_KwokTchelepi,hamon2016ordering,EfficientReorderedNonlinearGaussSeidelSolversWithHigherOrderForBlackOilModels_KlemetsdalRamussenMoynerLie,RobustNonlinearNewtonSolverWithAdaptiveInterfaceLocalizedTrustRegions_KlemetsdalMoynerLie} accelerate nonlinear convergence thanks to a reordering technique based on the flow direction.
The reordered systems have block-triangular structure and can therefore be efficiently solved with backward substitution.
To obtain highly scalable solvers able to take large time steps, multilevel solution algorithms relying on the application of the multigrid principles at the nonlinear level in a Full Approximation Scheme (FAS) have been proposed \cite{FullApproximationSchemeForReservoirSimulation_ToftLieMoyner,NonlinearMultigridSolversExploitingAmgeCoarseSpacesWithApproximationProperties_ChristensenVassilevskiVilla,NonlinearMultigridBasedOnLocalSpectralCoarseningForHeterogeneousDiffusionProblems_LeeHamonCastellettoVassulevskiWhite,AnAggregationBasedNonlinearMultigridSolverForTwoPhaseFlowAndTransportInPorousMedia_LeeHamonCastellettoVassulevskiWhite}.

The efforts to design efficient nonlinear preconditioners to the Newton update are particularly relevant to this work.
Two distinct directions have been explored. 
A class of nonlinear preconditioners leverages domain decomposition methods \cite{NonlinearAdditiveSchwarzPreconditionersAndApplicationInComputationalFluidDynamics_CaiKeyesMarcinkowski,HowToUse_aNonlinearSchwartzMethodToPreconditionNewtonMethod_DoleanMasson,AnAdaptiveSFDomainDecompositionSolver_KlemetsdalMoncorgeNilsenMoynerLie,NonlinearDomainDecompositionSchemeForSFIformulationOfCompositionalMultiphaseFlow_MoynerMoncorge,aNumericalStudyOfTheASPENasANonlinearpreconditionerForImmiscibleAndCompositionalPorousMediaFlow_KlemetsdalMoncorgeMoynerLie} to precondition the nonlinear system in a computationally inexpensive way and speed up convergence. 
In this work, we focus on nonlinear preconditioners obtained by splitting the system by physical field \cite{NonlinearPreconditioningForTwoPhaseFlows_LuoCaiKeyes,FullyImplicithybridTwoLevelDomainDecompositionAlgorithmsTorTwoPhaseFlowsInPorousMediaOn3DunstructuredGrids_LuoLiuCaiKeyes,NonlinearPreconditioningStrategiesForTwoPhaseFlowsInPorousMediaDiscretizedByAfullyImplicitDiscontinuousGalerkinMethod}.
Instead of decomposing the domain in space, field-split preconditioners exploit the mathematical structure of the nonlinear system to split the problem into multiple subproblems solved nonlinearly at a loose tolerance before the computation of a global update for all the degrees of freedom.
This approach provides an efficient framework to precondition the system field by field to compute more accurate nonlinear updates and in turn, accelerate nonlinear convergence \cite{Two-ScalePreconditioningForTwo-PhaseNonlinearFlowsInPorousMedia_SkogestadKeilegavlenNordbotten,SIN_Yang}.
In the context of mixed elliptic-hyperbolic multiphase flow and transport in porous media, the method is particularly attractive since it retains the best features of both SFI and FIM. 
Specifically, the preconditioning step resembles the SFI outer iteration based on a pressure update followed by a saturation update and therefore does not involve solving a large coupled linear system.
The global step of field-split algorithms plays the role of the computation of the Newton update in FIM and maintains robust convergence properties for strongly coupled problems.

In this work, we focus on immiscible two-phase flow in porous media and compare various field-split preconditioners based on the pressure-saturation decomposition. 
We consider the Multiplicative Schwarz Preconditioned for Inexact Newton (MSPIN) proposed in \cite{FieldSplitPreconditionedInexactNewtonAlgorithms_LiuKeyes}, in which the solution of the preconditioning step is used to compute an approximation of the preconditioned Jacobian system.
The global step consists in solving this global preconditioned system to obtain the updated pressures and saturations.
We also propose an alternative method, Field-Split Multiplicative Schwarz Newton (FSMSN) in which the preconditioning step is the same as in MSPIN, but the global step is simply the Newton update computed by linearizing the system at the preconditioned state.
We compare MSPIN and FSMSN to SFI and Newton with damping for FIM on challenging two-phase flow test cases. 
The comparison takes into account the role of the upwinding scheme---Phase Potential Upwinding \cite{AnAnalysisOfUpstreamDifferencing_Sammon,UpstreamDifferencingForMultiphaseFlowInReservoirSimulation_BrenierJaffre} or Implicit Hybrid Upwinding \cite{HybridUpwindDiscretizationOfNonlinearTwoPhaseFlowWithGravity_LeeEfendievTchelpi,HybridDiscretizationOfMultiPhaseFlowInPorousMediaInThePresenceOfViscousGravitationalAndCapillary_LeeEfendiev,ImplicitHybridUpwindSchemeForCoupledMultiphaseFlowAndTransportWithBuoyancy_HamonMallisonTchelepi,ConsistentUpwindingForSFImultiscaleCompositionalSimulation_MoncorgeMoynerTchelepi,AnalysisOfHybridUpwindingForFullyImplicitSimulationOfThreePhaseFlowWithGravity_HamonTchelepi,FullyImplicitMultidimensionalHybridUpwindSchemeForCoupledFlowAndTransport_HamonMallion,SmoothImplicitHybridUpwindingForCompositionalMultiphaseFlowInPorousMedia_BosmaHamonMallisonTchelepi,VertexApproximateGradientDiscretizationPreservingPositivityForTwoPhaseDarcyFlowsInIeterogeneousPorousMedia_BrennerMassonQuenjel}---as well as the role of the nonlinear tolerance used in the subproblems.
To do so, we propose an adaptive method to compute the subproblem nonlinear tolerance at each preconditioning step, as typically done in inexact Newton methods \cite{eisenstat96,DemboEisenstatSteihaug82, eisenstat94,InexactMethodsForBlackOilSequentialFullyImplicitSFIscheme_ZhouJiangTomin,InexactMethodsForSFIreservoirSimulation_JiangTominZhou,InexactNewtonMethodForGeneralPurposeReservoirSimulation_ShethMoncorge}.
We demonstrate that the two-step nonlinear preconditioners and, in particular, FSMSN, are successful at accelerating nonlinear convergence and are worth being explored as viable options to reduce the computational cost of reservoir simulation--which will be addressed in a future publication.

The structure of the present article is as follows.
We review the PDEs governing two-phase flow and transport in porous media in Section~\ref{sec:equations}. 
The fully implicit finite-volume scheme is reviewed in Section~\ref{sec:discretization_scheme}. 
The two nonlinear preconditioning techniques considered here--namely, FSMSN and MSPIN--are reviewed in Sections~\ref{sec:FSMSN} and~\ref{sec:preconditionedMethod},
respectively.
The method used to compute the adaptive subproblem tolerance is presented in Section~\ref{sec:adaptive_tolerance}.
We compare the efficiency of the nonlinear preconditioners with numerical examples in Section~\ref{sec:numericalResults}.

%%%%%%%%%%%%%%%%%%%%%%%%%%%%%%%%%%%%%%%%%%%%%%%%%%%%%%%%%%%%%%%%%%%%%%%%%%%%%%%%%%%%%%% 
%                           Equations
%%%%%%%%%%%%%%%%%%%%%%%%%%%%%%%%%%%%%%%%%%%%%%%%%%%%%%%%%%%%%%%%%%%%%%%%%%%%%%%%%%%%%%%

\section{Governing equations}\label{sec:equations}

Let $\overline{\Omega} = \Omega \cup \Gamma $ be a closed set in $\mathbb{R}^2$, with $\Omega$ an open set and $\Gamma$ its boundary.
We note $\mathbb{I} = [0,t_{max}]$ a finite time interval with $t_{max} > 0$ the maximal time.
We denote the spatial coordinates by $\bx \in \Omega$ and the time coordinate by $t\in \mathbb{I}$.
We consider the immiscible flow of two fluid phases in an incompressible porous medium, with a wetting ($w$) and a non-wetting ($nw$) phase.
For a phase $\ell \in \{ nw, w \}$, the mass conservation equation reads
\begin{equation}
  \label{eq:massConservation}
  \phi \dfrac{\partial ( \rho_\ell s_\ell)}{\partial t} 
  +\nabla \cdot \left(\rho_\ell \boldsymbol{u}_\ell(p,s) \right) = q_\ell, \qquad \ell \in \{ nw,w  \}, \quad \text{on} \, \, \Omega \times \mathbb{I},
\end{equation}
where $\phi=\phi(\bx)$ is the porosity of the medium, $q_\ell = q_\ell(\bx)$ is the source/sink term with the convention that $q_\ell > 0$ for injection, and $q_\ell<0$ for production.
The saturation $s_\ell=s_\ell(\bx,t)$ represents the fraction of the pore volume occupied by phase $\ell$, with the following constraint:
\begin{equation}
  \label{eq:satConstraint}
  \sum_\ell s_\ell = 1.
\end{equation}
We choose the wetting-phase saturation as the primary saturation unknown, denoted by $s := s_w$.
Using the saturation constraint \eqref{eq:satConstraint}, we write all the saturation-dependent properties as a function of $s$ only.
The phase velocity $\bu_{\ell}$ of a phase $\ell$ is given by Darcy's law
\begin{equation}
\label{eq:DarcyLaw}
  \bu_\ell( p, s ) = -k \lambda_\ell (\nabla p - \rho_\ell g\nabla d), \qquad \ell \in \{ nw,w \}, 
\end{equation}
where we have neglected capillary pressure.
In \eqref{eq:DarcyLaw}, $p(\bx,t)$ is the pressure, $k(\bx)$ is the scalar rock permeability, $\rho_\ell$ is the phase density, and $\lambda_\ell(s) = k_{r\ell}(s) / \mu_\ell$ is the phase mobility, defined as the phase relative permeability, $k_{r\ell}(s)$, divided by the phase viscosity $\mu_\ell$.
The gravitational acceleration is $g$ and the depth is $d$ (positive going downward).
Inserting the expression of the phase velocity given by Darcy's law \eqref{eq:DarcyLaw} into the mass conservation equation \eqref{eq:massConservation} gives rise to the following form of the two-phase flow and transport equation:
\begin{equation}
  \label{eq:TwoPhaseFlowStandardForm}
  \phi \dfrac{\partial ( \rho_\ell s_\ell)}{\partial t} - \nabla \cdot \left(k \rho_{\ell} \lambda_\ell (\nabla p - \rho_\ell g \nabla d) \right) = q_\ell, \qquad  \ell \in \{ nw,w \}, \quad \text{on} \, \Omega \times \mathbb{I}.
\end{equation}

In the sequential and field-split solution methods considered in this work, we employ a discretization scheme applied to a split form of the governing equations consisting of an elliptic pressure equation coupled with a hyperbolic transport problem.
The pressure subproblem is obtained by summing the mass balance equations \eqref{eq:TwoPhaseFlowStandardForm} over the two phases.
Using the saturation constraint \eqref{eq:satConstraint} and assuming that the two phases are incompressible, we obtain 
\begin{equation}
    \label{eq:pressureEquation}
    \nabla \cdot \bu_T(p,s) =\sum_\ell q_\ell, 
\end{equation}
where the total velocity $\bu_T$ is defined by $\bu_T(p,s) := \sum_{\ell} \bu_{\ell}(p,s)$.
The transport problem is obtained by eliminating the pressure variable in the flux term of \eqref{eq:TwoPhaseFlowStandardForm} to obtain the following fractional flow formulation as explained in \cite{ImplicitHybridUpwindSchemeForCoupledMultiphaseFlowAndTransportWithBuoyancy_HamonMallisonTchelepi,hamon2018implicit}. Precisely, we first write that
\begin{equation}
    \bu_T(p,s) := \sum_{\ell} \bu_{\ell}(p,s) = \lambda_T (-k \nabla p ) + ( \lambda_w \rho_w + \lambda_{nw} \rho_{nw} ) g\nabla d,
\end{equation}
which gives an expression of $-k \nabla p$ as a function $\bu_T$ and the gravity weights
\begin{equation}
    - k \nabla p = \frac{1}{\lambda_T} \bigg( \bu_T(p,s) - ( \lambda_w \rho_w + \lambda_{nw} \rho_{nw} ) g\nabla d \bigg).
\end{equation}
Using this in Darcy's law, we obtain
\begin{equation}
    \label{eq:fracFlowFormulation}
    \phi \dfrac{\partial (\rho_{\ell} s_{\ell})}{\partial t}+ \nabla \cdot \left( \rho_{\ell}\dfrac{\lambda_\ell}{\lambda_T} \bu_T(p,s) + k \rho_{\ell}\dfrac{\lambda_m\lambda_\ell}{\lambda_T} (\rho_\ell - \rho_m ) g \nabla d  \right) = q_\ell, \qquad m,\ell \in \{ nw, w \}, \, m \neq \ell,
\end{equation}
where $\lambda_T(s) := \sum_\ell \lambda_\ell(s)$ is the total mobility.

%%%%%%%%%%%%%%%%%%%%%%%%%%%%%%%%%%%%%%%%%%%%%%%%%%%%%%%%%%%%%%%%%%%%%%%%%%%%%%%%%%%%%
%                           Discretization
%%%%%%%%%%%%%%%%%%%%%%%%%%%%%%%%%%%%%%%%%%%%%%%%%%%%%%%%%%%%%%%%%%%%%%%%%%%%%%%%%%%%%

\section{Discretization scheme}\label{sec:discretization_scheme}

Given a mesh consisting of $M$ cells, we use first-order finite-volume scheme.
We consider $0= t^0 \leq t^1 \leq \cdots \leq t^N = t_{max}, \,N \in \mathbb{N}$ a finite discretization of the temporal axis, and $\Delta t^{n+1} := t^{n+1}-t^n, \, 0\leq n \leq N, \, n\in \mathbb{N}$ the time step.
We use a (backward-Euler) fully implicit scheme for the integration in time.
We denote by $\mathbf{p}^{n+1} := [p^{n+1}_1, \dots, p^{n+1}_M]^{\intercal}$ and $\mathbf{s}^{n+1} := [ s^{n+1}_1, \dots, s^{n+1}_M]^{\intercal}$, the solution pair $( \mathbf{p}^{n+1}, \mathbf{s}^{n+1} )$ the vectors collecting respectively the pressure and saturation unknowns.
To define the discrete problem, we first introduce the phase-based residual in cell $K$ at time $n+1$, $r^{n+1}_{\ell,K}$, as
\begin{equation}
  r^{n+1}_{\ell,K}( \mathbf{p}^{n+1}, \mathbf{s}^{n+1} ) := V_K \phi_K \dfrac{ \rho^{n+1}_{\ell,K} s_{\ell,K}^{n+1}-\rho^n_{\ell,K} s_{\ell,K}^n}{\Delta t^{n+1}} + \sum_{L\in adj(K)} F_{\ell,KL}^{n+1} ( \mathbf{p}^{n+1}, \mathbf{s}^{n+1} ) - V_K q_{\ell,K}( p^{n+1}_K, s^{n+1}_K ).
  \label{eq:definition_discrete_phase_based_residual}
\end{equation}
In \eqref{eq:definition_discrete_phase_based_residual}, $F_{\ell,KL}^{n+1}$ is the numerical flux for the interface $(KL)$ between cells $K$ and $L$, $\textit{adj}(K)$ is the set of neighbors of cell $K$ and $V_K$ the volume of cell $K$.
The computation of the numerical flux is based on a Two-Point Flux Approximation (TPFA).
We consider both Phase-Potential Upwinding (PPU) \cite{AnAnalysisOfUpstreamDifferencing_Sammon,UpstreamDifferencingForMultiphaseFlowInReservoirSimulation_BrenierJaffre} and Implicit Hybrid Upwinding (IHU) \cite{HybridUpwindDiscretizationOfNonlinearTwoPhaseFlowWithGravity_LeeEfendievTchelpi,HybridDiscretizationOfMultiPhaseFlowInPorousMediaInThePresenceOfViscousGravitationalAndCapillary_LeeEfendiev,ImplicitHybridUpwindSchemeForCoupledMultiphaseFlowAndTransportWithBuoyancy_HamonMallisonTchelepi,ConsistentUpwindingForSFImultiscaleCompositionalSimulation_MoncorgeMoynerTchelepi,AnalysisOfHybridUpwindingForFullyImplicitSimulationOfThreePhaseFlowWithGravity_HamonTchelepi,FullyImplicitMultidimensionalHybridUpwindSchemeForCoupledFlowAndTransport_HamonMallion,alali2021finite,SmoothImplicitHybridUpwindingForCompositionalMultiphaseFlowInPorousMedia_BosmaHamonMallisonTchelepi} for the approximation of the transport coefficient of the flux term.

In this work, we consider two equivalent formulations of the nonlinear problem.
In the Newton-based FIM, we apply the standard fully coupled Newton's method with damping to find the solution pair $( \mathbf{p}^{n+1}, \mathbf{s}^{n+1} )$ that satisfies
\begin{equation}
  \mathbf{r}^{n+1}( \mathbf{p}^{n+1}, \mathbf{s}^{n+1} ) := [ \mathbf{r}^{n+1}_{nw}( \mathbf{p}^{n+1}, \mathbf{s}^{n+1} ), \mathbf{r}^{n+1}_w( \mathbf{p}^{n+1}, \mathbf{s}^{n+1} ) ]^{\intercal} = \mathbf{0}.
    \label{eq:definition_full_phase_based_residual}
\end{equation}
In sequential and field-split nonlinear solution strategies, we split the problem into a discrete pressure problem and a discrete transport problem.
The solution pair $( \mathbf{p}^{n+1}, \mathbf{s}^{n+1} )$ must in this case satisfy:
\begin{equation}
  \mathbf{f}^{n+1}( \mathbf{p}^{n+1}, \mathbf{s}^{n+1} ) := [ \mathbf{g}^{n+1}( \mathbf{p}^{n+1}, \mathbf{s}^{n+1} ), \mathbf{h}^{n+1}( \mathbf{p}^{n+1}, \mathbf{s}^{n+1} ) ]^{\intercal} = \mathbf{0},
  \label{eq:definition_full_field_split_residual}
\end{equation}
where the discrete pressure problem is obtained by summing the phase-based residual equations, as in \eqref{eq:pressureEquation}:
\begin{equation}
  \mathbf{g}^{n+1}( \mathbf{p}^{n+1}, \mathbf{s}^{n+1} ) := \mathbf{r}^{n+1}_{nw}( \mathbf{p}^{n+1}, \mathbf{s}^{n+1} ) + \mathbf{r}^{n+1}_w( \mathbf{p}^{n+1}, \mathbf{s}^{n+1} ),
\end{equation}
and the discrete saturation problem simply consists of one of the phase-based residual equations. Here, we choose the wetting-phase residual:
\begin{equation}
  \mathbf{h}^{n+1}( \mathbf{p}^{n+1}, \mathbf{s}^{n+1} ) := \mathbf{r}^{n+1}_w( \mathbf{p}^{n+1}, \mathbf{s}^{n+1} ).
\end{equation}
From now on, we consider time step $n \in \mathbb{N}$ and describe two algorithms to compute the solution pair $( \mathbf{p}^{n+1}, \mathbf{s}^{n+1} )$ of nonlinear problems \eqref{eq:definition_full_phase_based_residual} and \eqref{eq:definition_full_field_split_residual} at time $n+1$.
For simplicity, we drop the temporal superscript denoting the time step.

%%%%%%%%%%%%%%%%%%%%%%%%%%%%%%%%%%%%%%%%%%%%%%%%%%%%%%%%%%%%%%%%%%%%%%%%%%%%%%%%%%%%%
%                           FSMSN
%%%%%%%%%%%%%%%%%%%%%%%%%%%%%%%%%%%%%%%%%%%%%%%%%%%%%%%%%%%%%%%%%%%%%%%%%%%%%%%%%%%%%

\section{Field-split multiplicative Schwarz Newton method}
\label{sec:FSMSN}
The multiplicative Schwarz method has been used to split boundary value problems (BVP) into subproblems solver on smaller physical domains \cite{NonlinearAdditiveSchwarzPreconditionersAndApplicationInComputationalFluidDynamics_CaiKeyesMarcinkowski,HowToUse_aNonlinearSchwartzMethodToPreconditionNewtonMethod_DoleanMasson,AnAdaptiveSFDomainDecompositionSolver_KlemetsdalMoncorgeNilsenMoynerLie,NonlinearDomainDecompositionSchemeForSFIformulationOfCompositionalMultiphaseFlow_MoynerMoncorge,aNumericalStudyOfTheASPENasANonlinearpreconditionerForImmiscibleAndCompositionalPorousMediaFlow_KlemetsdalMoncorgeMoynerLie}.
It has also been used to split a coupled BVP into subproblems based on the physics \cite{NonlinearPreconditioningForTwoPhaseFlows_LuoCaiKeyes,FullyImplicithybridTwoLevelDomainDecompositionAlgorithmsTorTwoPhaseFlowsInPorousMediaOn3DunstructuredGrids_LuoLiuCaiKeyes,NonlinearPreconditioningStrategiesForTwoPhaseFlowsInPorousMediaDiscretizedByAfullyImplicitDiscontinuousGalerkinMethod}, each subproblem being solved on the full domain to update one of the fields (here, pressure and saturation).
The motivation for splitting a coupled problem according to the physics is to solve the physical subproblems one at a time (ideally, with a specialized solver) and use the individually updated fields to precondition the nonlinear iteration, yielding a faster nonlinear convergence.
In this section, the objective is to use the field-split approach to construct a predictor-corrector method that converges faster than commonly used algorithms based on the Newton iteration or the sequential fully implicit iteration.

The preconditioning step of the FSMSN outer iteration consists in computing an intermediate value of the pressure and saturation, $\mathbf{p}^{k,\star}$ and $\mathbf{s}^{k,\star}$, by solving individual pressure and transport problems sequentially.
We first solve nonlinearly the pressure equation.
Specifically, for a fixed $\mathbf{s}^k$, find the update $\boldsymbol{\delta}_1( \mathbf{p}^{k}, \mathbf{s}^{k} )$ such that %$\mathbf{p}^{k,\star}$ such that
\begin{equation}
  \mathbf{g}( \mathbf{p}^{k,\star}, \mathbf{s}^k ) = \mathbf{0},
  \label{eq:pressure_problem}
\end{equation}
where $\mathbf{p}^{k,\star} := \mathbf{p}^{k} + \boldsymbol{\delta}_1( \mathbf{p}^{k}, \mathbf{s}^{k} )$.
We solve \eqref{eq:pressure_problem} with Newton's method to obtain a pressure prediction, $\mathbf{p}^{k,\star}$.
The preconditioning step continues by solving nonlinearly a transport problem.
For a fixed $\mathbf{p}^{k,\star}$, find $\boldsymbol{\delta}_2( \mathbf{p}^{k},  \mathbf{s}^{k} )$ such that
\begin{equation}
  \mathbf{h}( \mathbf{p}^{k,\star}, \mathbf{s}^{k,\star} ) = \mathbf{0},
  \label{eq:transport_problem}
\end{equation}
where $\mathbf{s}^{k,\star} := \mathbf{s}^{k} + \boldsymbol{\delta}_2( \mathbf{p}^{k}, \mathbf{s}^{k} )$.
In \eqref{eq:transport_problem}, we explore two formulations of the discrete transport problem.
In the first formulation, we use directly the intermediate pressure, $\mathbf{p}^{k,\star}$, and we consider a discrete transport equation approximating \eqref{eq:TwoPhaseFlowStandardForm} for $\ell = w$.
In the second formulation, we compute an intermediate total velocity field with the intermediate pressure, and we consider a discrete transport equation approximation \eqref{eq:fracFlowFormulation} for $\ell = w$ with a fixed total velocity.
This technique is commonly used in the sequential fully implicit method in which the total velocity is also fixed during the resolution of the transport problem.
In both options, we use Newton's method with damping to solve the discrete transport problem, although more efficient nonlinear solvers are available \cite{natvig2008fast,PotentialBasedReducedNewtonAlgorithmForNonlinearMultiphaseFlowInPorousMedia_KwokTchelepi,hamon2016ordering,EfficientReorderedNonlinearGaussSeidelSolversWithHigherOrderForBlackOilModels_KlemetsdalRamussenMoynerLie,RobustNonlinearNewtonSolverWithAdaptiveInterfaceLocalizedTrustRegions_KlemetsdalMoynerLie}.

In the correction step, we re-evaluate the residual and Jacobian matrix with the intermediate pressure and saturation.
Then, we compute the $(k+1)$-th solution iterate as
\begin{equation}
\begin{aligned}
  \mathbf{p}^{k+1} &\leftarrow \mathbf{p}^{k,\star} + \boldsymbol{\delta} \mathbf{p}^{k+1}, \\
  \mathbf{s}^{k+1} &\leftarrow \mathbf{s}^{k,\star} + \boldsymbol{\tau}^{k+1} \boldsymbol{\delta} \mathbf{s}^{k+1}, 
\end{aligned}
\label{eq:fsmsn_update_application}
\end{equation}
where $\boldsymbol{\tau}^{k+1}$ is a diagonal matrix of damping parameters.
In \eqref{eq:fsmsn_update_application}, the update is obtained by solving the linear system
\begin{equation}
  \boldsymbol{\delta} \mathbf{x}^{k+1} = - \mathbf{J}( \mathbf{x}^{k,\star} )^{-1} \mathbf{r}( \mathbf{x}^{k,\star} ), \qquad \mathbf{x} := \begin{pmatrix} \mathbf{p}\\ \mathbf{s} \end{pmatrix}.
  \label{eq:fsmsn_global_solve}
\end{equation}
Note that for the assembly of $\mathbf{J}$ in \eqref{eq:fsmsn_global_solve}, the fluxes are computed using the same discretization as in the subproblems \eqref{eq:pressure_problem} and \eqref{eq:transport_problem} to maintain a uniform discrete approach throughout the FSMSN nonlinear iteration.
The implementation of FSMSN is described in Algorithm~\ref{algo:FSMSN}.  
%
%In summary, at each iteration of FSMSN, we first predict the pressure and saturation update with the SFI method, and then correct the coupled results with Newton's method. %As a combination of two schemes that converge to the same solution for a well-chosen initial guess, we conclude that the FSMSN method converges to the same solution.

\begin{algorithm}
  \caption{FSMSN algorithm for two-phase flow and transport}
  \label{algo:FSMSN}
  \begin{algorithmic}
    \FOR{$k = 1, \dots, k_{\textit{itermax}}$}
    \STATE{Check convergence, and break nonlinear loop if convergence was achieved}
    \STATE{\textit{Pressure step:}}
    \STATE{$\quad$ For a fixed $\mathbf{s}^k$, solve $\mathbf{g}( \mathbf{p}^{k,\star}, \mathbf{s}^k ) = \mathbf{0}$}
    \STATE{\textit{Transport step:}} 
    \STATE{$\quad$ For a fixed $\mathbf{p}^{k,\star}$ or a fixed total velocity $\mathbf{u}^{k,\star}_T$, solve $\mathbf{h}( \mathbf{p}^{k,\star}, \mathbf{s}^{k,\star} ) = \mathbf{0}$}
    \STATE{\textit{Coupled step:}}
    \STATE{$\quad$ Recompute residual and Jacobian using $(\mathbf{p}^{k,\star}, \mathbf{s}^{k,\star})$}    
    \STATE{$\quad$ Solve linear system: $\boldsymbol{\delta} \mathbf{x}^{k+1} = - \mathbf{J}( \mathbf{x}^{k,\star} )^{-1} \mathbf{r}( \mathbf{x}^{k,\star} )$}
    \STATE{$\quad$ Update solution: $\mathbf{p}^{k+1} \leftarrow \mathbf{p}^{k,\star} + \boldsymbol{\delta} \mathbf{p}^{k+1}$; $\mathbf{s}^{k+1} \leftarrow \mathbf{s}^{k,\star} + \tau^{k+1} \boldsymbol{\delta} \mathbf{s}^{k+1}$}
    \ENDFOR
  \end{algorithmic}
\end{algorithm}

%%%%%%%%%%%%%%%%%%%%%%%%%%%%%%%%%%%%%%%%%%%%%%%%%%%%%%%%%%%%%%%%%%%%%%%%%%%%%%%%%%%%%
%                           Field split Preconditioning methods
%%%%%%%%%%%%%%%%%%%%%%%%%%%%%%%%%%%%%%%%%%%%%%%%%%%%%%%%%%%%%%%%%%%%%%%%%%%%%%%%%%%%%
\section{Preconditioned method based on field-split multiplicative Schwarz method}
\label{sec:preconditionedMethod}

In this section, we review the construction of a preconditioned method for the two-phase problem using a field-split multiplicative Schwarz method (MSPIN).
As in the previous section, the iteration consists of two steps, with first a preconditioning step following with a global step.
The preconditioning step is identical to that of FSMSN.
Using the same notations as in Section~\ref{sec:FSMSN}, the flow problem consists in finding a new pressure field, $\mathbf{p}^{k,\star} := \mathbf{p}^{k} + \boldsymbol{\delta}_1( \mathbf{p}^{k}, \mathbf{s}^{k} )$, that satisfies \eqref{eq:pressure_problem} for a fixed saturation field, $\mathbf{s}^{k}$.
Then, the transport problem consists in finding an updated saturation field, $\mathbf{s}^{k,\star} := \mathbf{s}^{k} + \boldsymbol{\delta}_2( \mathbf{p}^{k}, \mathbf{s}^{k} )$, that satisfies \eqref{eq:transport_problem} for a fixed pressure field, $\mathbf{p}^{k,\star}$.
In this second step, we only consider the formulation in which the discrete transport equation approximates \eqref{eq:TwoPhaseFlowStandardForm} for $\ell = w$.
The MSPIN algorithm differs from FSMSN in the global step.
Once the pressure and transport problems are solved, we form the following preconditioned problem:
\begin{equation}
  \label{eq:MSPIN_problem}
  \boldsymbol{\set{F}}(\mathbf{p}^k,\mathbf{s}^k) := [\boldsymbol{\delta}_1(\mathbf{p}^k,\mathbf{s}^k), \boldsymbol{\delta}_2(\mathbf{p}^k,\mathbf{s}^k)]^{\intercal} = \mathbf{0}.
\end{equation}
The goal of the coupled step is to form a Jacobian system from \eqref{eq:MSPIN_problem} and solve it to obtain the pressure and saturation updates of outer iteration $k$.
Forming the Jacobian system requires computing the partial derivatives of $\boldsymbol{\delta}_1$ and $\boldsymbol{\delta}_2$ with respect to the pressure and saturation variables.
Using the chain rule, we obtain the derivatives of $\boldsymbol{\delta}_1$ by differentiating \eqref{eq:pressure_problem}, which yields:
\begin{align}
  &\partial_p \boldsymbol{\delta}_1(\mathbf{p}^k,\mathbf{s}^k) = - \mathbf{I}_M, \label{eq:ddelta_1_dp} \\ 
  &\partial_s \boldsymbol{\delta}_1(\mathbf{p}^k,\mathbf{s}^k) = - \big(  \partial_p \mathbf{g}(\mathbf{p}^{k,\star},\mathbf{s}^k) \big)^{-1} \partial_s \mathbf{g}(\mathbf{p}^{k,\star},\mathbf{s}^k). \qquad \qquad \qquad  \qquad \qquad \qquad \qquad \quad \label{eq:ddelta_1_ds}
\end{align}
$\mathbf{I}_M$ denotes a $M$ by $M$ identity matrix. Differentiating \eqref{eq:transport_problem} yields the derivatives of $\boldsymbol{\delta}_2$:
\begin{align}
  &\partial_p \boldsymbol{\delta}_2(\mathbf{p}^k,\mathbf{s}^k) = \mathbf{0}, \label{eq:ddelta_2_dp} \\
  &\partial_s \boldsymbol{\delta}_2(\mathbf{p}^k,\mathbf{s}^k) = - \mathbf{I}_M - \big( \partial_s \mathbf{h}(\mathbf{p}^{k,\star},\mathbf{s}^{k,\star}) \big)^{-1} \partial_p \mathbf{h}(\mathbf{p}^{k,\star},\mathbf{s}^{k,\star}) \big( \partial_p \mathbf{g}(\mathbf{p}^{k,\star},\mathbf{s}^k) \big)^{-1} \partial_s \mathbf{g}(\mathbf{p}^{k,\star},\mathbf{s}^k). \label{eq:ddelta_2_ds}
\end{align}
Using equations \eqref{eq:ddelta_1_dp} to \eqref{eq:ddelta_2_ds}, the Jacobian matrix, $\boldsymbol{\set{J}}(\mathbf{p}^k,\mathbf{s}^k)$ of \eqref{eq:MSPIN_problem} is then given by:
\begin{equation}
  \label{eq:Jacobian_prec_func}
  \begin{aligned}
    \boldsymbol{\set{J}}(\mathbf{p}^k,\mathbf{s}^k) & = -
    \begin{pmatrix}
      \partial_p \mathbf{g}(\mathbf{p}^{k,\star},\mathbf{s}^k) & \mathbf{0}\\
      \partial_p \mathbf{h}(\mathbf{p}^{k,\star},\mathbf{s}^{k,\star}) & \partial_s \mathbf{h}(\mathbf{p}^{k,\star},\mathbf{s}^{k,\star})
    \end{pmatrix}^{-1}
    \begin{pmatrix}
      \partial_p \mathbf{g}(\mathbf{p}^{k,\star},\mathbf{s}^k) & \partial_s \mathbf{g}(\mathbf{p}^{k,\star},\mathbf{s}^k) \\
      \partial_p \mathbf{h}(\mathbf{p}^{k,\star},\mathbf{s}^{k,\star}) & \partial_s \mathbf{h}(\mathbf{p}^{k,\star},\mathbf{s}^{k,\star}).
    \end{pmatrix}
  \end{aligned}
\end{equation}
At this point, the second matrix on the right-hand side of \eqref{eq:Jacobian_prec_func} differs from the unpreconditioned Jacobian matrix of \eqref{eq:definition_full_phase_based_residual} because the first row is evaluated at $(\mathbf{p}^{k,\star},\mathbf{s}^k)$, while the second row is evaluated at $(\mathbf{p}^{k,\star},\mathbf{s}^{k,\star})$.
Following the methodology of \cite{FieldSplitPreconditionedInexactNewtonAlgorithms_LiuKeyes}, we approximate the preconditioned Jacobian matrix of the Multiplicative Schwarz Preconditioned Inexact Newton (MSPIN) method by setting $\boldsymbol{\delta}_2( \mathbf{p}^k, \mathbf{s}^k ) = \mathbf{0}$ in \eqref{eq:Jacobian_prec_func}, which results in:
\begin{equation}
  \label{eq:MSPIN_approx_Jacobian}
  \boldsymbol{\set{J}}(\mathbf{p}^k,\mathbf{s}^k) \approx - \begin{pmatrix}
    \partial_p \mathbf{g}(\mathbf{p}^{k,\star},\mathbf{s}^k) & \mathbf{0}\\
    \partial_p \mathbf{h}(\mathbf{p}^{k,\star},\mathbf{s}^{k}) & \partial_s \mathbf{h}(\mathbf{p}^{k,\star},\mathbf{s}^k)
  \end{pmatrix}^{-1}
  \mathbf{J}(\mathbf{p}^{k,\star},\mathbf{s}^k).
\end{equation}
where $\mathbf{J}$ is the original, unpreconditioned Jacobian matrix of \eqref{eq:definition_full_phase_based_residual}.
We note that in \eqref{eq:MSPIN_approx_Jacobian}, the two matrices on the right-hand side are now fully evaluated using the state of the system after the pressure step.
Then, we obtain the $(k+1)$-th solution iterate as
\begin{equation}
\begin{aligned}
  \mathbf{p}^{k+1} &\leftarrow \mathbf{p}^{k} + \boldsymbol{\delta} \mathbf{p}^{k+1}, \\
  \mathbf{s}^{k+1} &\leftarrow \mathbf{s}^{k} + \boldsymbol{\tau}^{k+1} \boldsymbol{\delta} \mathbf{s}^{k+1}, 
\end{aligned}
\label{eq:mspin_update_application}
\end{equation}
In \eqref{eq:mspin_update_application}, the update is obtained by solving the linear system:
\begin{equation}
\boldsymbol{\delta} \mathbf{x}^{k+1} =  -\boldsymbol{\set{J}}(\mathbf{p}^k,\mathbf{s}^k)^{-1} \boldsymbol{\set{F}}(\mathbf{p}^k,\mathbf{s}^k), \qquad 
\mathbf{x} := \begin{pmatrix} \mathbf{p}\\ \mathbf{s} \end{pmatrix}.
  \label{eq:mspin_global_solve}
\end{equation}
The MSPIN iteration for two-phase flow and transport is summarized in Algorithm~\ref{algo:mspin}.
The MSPIN strategy converges to the same solution as Newton based FIM, as proven in \cite{FieldSplitPreconditionedInexactNewtonAlgorithms_LiuKeyes}.
In \cite{FieldSplitPreconditionedInexactNewtonAlgorithms_LiuKeyes}, the authors provide numerical tests based on the Navier-Stokes equations showing that MSPIN is more robust than Newton's method and the Additive Schwarz-type preconditioning (ASPIN).
In this work, we focus on the assessment of the nonlinear behavior of the MSPIN algorithm for the strongly coupled and highly nonlinear two-phase flow and transport problem.
In our implementation, we form matrix \eqref{eq:MSPIN_approx_Jacobian} and solve \eqref{eq:mspin_global_solve} by calling a direct solver twice.
In future work, we will exploit the structure of the block-triangular matrix of \eqref{eq:MSPIN_approx_Jacobian} inside an iterative Krylov-type linear solver \cite{GMRES_Saad_1986,HybridKrylovMethods_BrwonSaad_1990} to improve the efficiency of the approach on large-scale problems.

%In this section, we show how to construct a preconditioned method for the two-phase flow and transport problem using a field-split multiplicative Schwarz method. To this end, for each outer-loop iteration $k$, we recall the discrete two-phase flow pressure and transport problem:

%%%%%%%%%%%%%%%%%%%%%%%%%%%%%%%%%%%%%%%%%%%%%%%%%%%%%%%%%%%%%%%%%%%%%%%%%%%%%%%%%%%%%
%%
 
\begin{algorithm}
  \caption{MSPIN algorithm for two-phase flow and transport}
  \label{algo:mspin}
  \begin{algorithmic}
    \FOR{$k=1,\dots,k_{itermax}$}
    \STATE{Check convergence, and break nonlinear loop if convergence was achieved}
    \STATE{\textit{Pressure step:}}
    \STATE{$\quad$ For a fixed $\mathbf{s}^k$, solve $\mathbf{g}(\mathbf{p}^{k,\star},\mathbf{s}^k) = \mathbf{0}$}
    \STATE{$\quad$ Update and save $\partial_p \mathbf{g}(\mathbf{p}^{k,\star},\mathbf{s}^k)$, $\partial_s \mathbf{g}(\mathbf{p}^{k,\star},\mathbf{s}^k)$, $\partial_p \mathbf{h}(\mathbf{p}^{k,\star},\mathbf{s}^k)$, and $\partial_s \mathbf{h}(\mathbf{p}^{k,\star},\mathbf{s}^k)$} 
    \STATE{\textit{Transport step:}}
    \STATE{$\quad$ For fixed $\mathbf{p}^{k,\star}$, solve $ \mathbf{h}(\mathbf{p}^{k,\star},\mathbf{s}^{k,\star}) = \mathbf{0}$}
    \STATE{\textit{Coupled step:}}
    \STATE{$\quad$ Form preconditioned residual $\boldsymbol{\set{F}}(\mathbf{p}^k,\mathbf{s}^k)$} and Jacobian matrix $\boldsymbol{\set{J}}(\mathbf{p}^k,\mathbf{s}^k)$ as in \eqref{eq:MSPIN_problem} and \eqref{eq:Jacobian_prec_func}
    \STATE{$\quad$ Solve linear system: $\boldsymbol{\delta} \mathbf{x}^{k+1} =  -\boldsymbol{\set{J}}(\mathbf{p}^k,\mathbf{s}^k)^{-1} \boldsymbol{\set{F}}(\mathbf{p}^k,\mathbf{s}^k)$}
    \STATE{$\quad$ Update solution: $\mathbf{p}^{k+1} \leftarrow \mathbf{p}^k + \boldsymbol{\delta} \mathbf{p}^{k+1}$; $\mathbf{s}^{k+1}\leftarrow \mathbf{p}^k + \boldsymbol{\tau}^{k+1} \boldsymbol{\delta} \mathbf{s}^k$}
    \ENDFOR
  \end{algorithmic}
\end{algorithm}

\section{Convergence check and adaptive nonlinear tolerance}
\label{sec:adaptive_tolerance}

In Algorithms~\ref{algo:FSMSN} and \ref{algo:mspin}, the convergence checks involved in the outer and inner loops are performed using the $\ell_2$-norm of the normalized residual.
Specifically, convergence of the full problem is achieved when:
\begin{equation}
  \max_{\ell} \left( \, \left|\left| \mathbf{diag} \big( \mathbf{r}^{k+1}_{\ell}( \mathbf{p}^{k+1}, \mathbf{s}^{k+1} ) \big) 
  \mathbf{diag} \big(
  \mathbf{m}_{\ell}( \mathbf{p}^{n} ) \big)^{-1} \right|\right|_{2,2} \, \right)  < \epsilon.
\label{eq:convergence_full_problem}  
\end{equation}
where we have used the normalizer $\mathbf{m}_{\ell}( \mathbf{p}^{n} ) = [ \rho^{n}_{\ell,1} V_1 \phi_1, \dots, \rho^{n}_{\ell,M} V_M \phi_M ]^{\intercal}$ evaluated at the previous converged time step, and where $\mathbf{diag}(\mathbf{x})$ is the diagonal matrix obtained from the argument vector $\mathbf{x}$.
In the first inner loop of Algorithms~\ref{algo:FSMSN} and \ref{algo:mspin}, convergence of the pressure problem is reached when 
\begin{equation}
  \left|\left| \mathbf{diag} \big(  \mathbf{g}( \mathbf{p}^{k,\star}, \mathbf{s}^k ) \big) \mathbf{diag} \big( \mathbf{m}_{\textit{nw}}( \mathbf{p}^{n} ) + \mathbf{m}_{\textit{w}}( \mathbf{p}^{n} ) \big)^{-1}  \right|\right|_{2,2}  < \epsilon^k_{p},
\label{eq:convergence_pressure_problem}  
\end{equation}
while, in the second inner loop, the transport problem is converged when
\begin{equation}
  \left|\left| \mathbf{diag} \big( \mathbf{h}( \mathbf{p}^{k,\star}, \mathbf{s}^{k,\star} ) \big) \mathbf{diag} \big( \mathbf{m}_{\textit{w}}( \mathbf{p}^{n} ) \big)^{-1} \right|\right|_{2,2}  < \epsilon^k_s.
\label{eq:convergence_transport_problem}  
\end{equation}
If these criteria are not satisfied, we keep iterating until we reach convergence or the maximum number of iterations.
The choice of subproblem tolerances in \eqref{eq:convergence_pressure_problem} and \eqref{eq:convergence_transport_problem} is key for the nonlinear behavior and the efficiency of the schemes.
In the numerical examples, we explore the two approaches detailed below.

In the first approach, we set the subproblem tolerances, $\epsilon^k_p$ and $\epsilon^k_s$, to constant values, $\epsilon_p$ and $\epsilon_s$, chosen to be stricter or equal to the full problem tolerance, $\epsilon$, used in the outer loop.
This is motivated by the fact that, in MSPIN, the computation of the preconditioned Jacobian in \eqref{eq:Jacobian_prec_func} assumes that the pressure and transport subproblems are fully converged and that \eqref{eq:pressure_problem} and \eqref{eq:transport_problem} are satisfied.
However, this approach requires a significant computational effort to solve the subproblems which is likely to undermine the efficiency of the scheme.

In an alternative approach, we also explore the use of adaptive subproblem tolerances to minimize the number of subproblem iterations and reduce the computational cost of the schemes.
We define the subproblem tolerances as:
\begin{align}
  &\epsilon^0_{p} = \epsilon^0_{s} := 1 \\ 
  &\epsilon^k_{p} := \eta^k \epsilon^{k-1}_{p}
  \qquad \text{and} \qquad
  \epsilon^k_{s} := \eta^k \epsilon^{k-1}_{s}, \quad k \geq 1.
\end{align}
The parameter $\eta^k \in [0,1[$ depends on the outer iteration number and is used to control the subproblem tolerance.
We choose a relatively large $\eta^k$ during the first outer iterations to use a relaxed tolerance, and we gradually reduce $\eta^k$ to obtain a tighter tolerance as the outer loop approaches convergence.
A similar approach based on the Eisenstat-Walker algorithm is commonly used to reduce the cost of solving the linear systems in the inexact Newton method (see for instance \cite{DemboEisenstatSteihaug82, eisenstat94}).
Computing the parameter $\eta^k$ at each outer iteration is the critical part of the algorithm.
We adapt the work of \cite{DemboEisenstatSteihaug82, eisenstat94} to Algorithms~\ref{algo:FSMSN} and \ref{algo:mspin} by setting:
\begin{align}
  &(A_1) \qquad \eta^k = 0.1, \quad k \geq 0, \label{eq:forcingTerm_A1}\\
  &(A_2) \qquad \eta^k = 2^{-(k+1)}, \quad k \geq 0, \label{eq:forcingTerm_A2}\\
  &(A_3) \qquad \eta^0 = 1, \qquad
   \eta^k = \dfrac{ \left| \, \| \mathbf{r}^{k}( \mathbf{p}^k, \mathbf{s}^k ) \|_2 - \| \mathbf{r}^{k-1}( \mathbf{p}^{k-1}, \mathbf{s}^{k-1} )  + \mathbf{J}( \mathbf{x}^{k-1} ) \delta \mathbf{x}^{k-1} \|_2 \, \right|}{\| \mathbf{r}^{k-1}( \mathbf{p}^{k-1}, \mathbf{s}^{k-1} ) \|_2}, \quad k \geq 1. \label{eq:forcingTerm_A3}
\end{align}
In the next section, we refer to these approaches as $A_1$, $A_2$, and $A_3$, respectively.
The algorithm employed to compute the subproblem tolerance and use in the inner loop is illustrated in Algorithm~\ref{algo:inexactNewton} for the transport problem.
\begin{algorithm}
  \caption{Transport inner loop with adaptive tolerance}
  \label{algo:inexactNewton}
  \begin{algorithmic}
    \STATE{Given $k \in \mathbb{N}^*$ the outer-loop iteration number}
    \STATE{Given $\epsilon^{k-1}_s \in [0,1[$ the tolerance used at outer-loop iteration $k-1$ ($\epsilon^0 = 1$)}
    \STATE{Compute $\eta^{k}$ using \eqref{eq:forcingTerm_A1}, or \eqref{eq:forcingTerm_A2}, or \eqref{eq:forcingTerm_A3}}
    \STATE{Compute the new tolerance $\epsilon^k_s = \eta^{k} \epsilon^{k-1}_s$}
    \FOR{$m = 1,\cdots,m_{maxiter}$}
      \STATE{Compute the residual, $\mathbf{h}( \mathbf{p}^{k,\star}, \mathbf{s}^{k,m} )$, and the Jacobian matrix}
      \IF{$|| \mathbf{h}( \mathbf{p}^{k,\star}, \mathbf{s}^{k,\star} ) / \mathbf{m}_{\textit{w}}( \mathbf{p}^{n} ) ||_{2}  < \epsilon^k_s$}
         \STATE{Convergence is achieved, return the solution $\mathbf{s}^{k,\star}$}
      \ENDIF
      \STATE{Solve the Jacobian system to compute $\delta \mathbf{s}^{k,m+1}$}
      \STATE{Update the saturation solution: $\mathbf{s}^{k,m+1} \leftarrow \mathbf{s}^{k,m} + \tau^{m+1} \delta \mathbf{s}^{k,m+1}$}
    \ENDFOR
  \end{algorithmic}
\end{algorithm}

%%%%%%%%%%%%%%%%%%%%%%%%%%%%%%%%%%%%%%%%%%%%%%%%%%%%%%%%%%%%%%%%%%%%%%%%%%%%%%%%%%%%
%                           Results
%%%%%%%%%%%%%%%%%%%%%%%%%%%%%%%%%%%%%%%%%%%%%%%%%%%%%%%%%%%%%%%%%%%%%%%%%%%%%%%%%%%%

\section{Numerical results}
\label{sec:numericalResults}

To compare the different algorithms discussed above, we consider various numerical examples with heterogeneous permeability fields.
For a given time step $\Delta t$, the local phase-based CFL number is defined as
\begin{equation}
  \label{eq:cfl}
  \text{CFL}_{\ell,K} := \dfrac{\Delta t \sum_{L\in adj(K)} \max \big( 0, F_{\ell,KL}^{n+1} (\mathbf{p}^{n+1}, \mathbf{s}^{n+1}) \big)}{V_K\phi_K},  \, \ell = {nw, w},
\end{equation}
where $V_K$ is the volume of the grid cell $K$, $\phi_K$ is the porosity and $F_{\ell,KL}$ is the numerical flux of the phase $\ell$ for the interface $(KL)$ between cells $K$ and $L$.
The maximum CFL number is then computed as
\begin{equation}
\text{max CFL} = \argmax_{\ell, 1\leq K\leq M} \text{CFL}_{\ell, K},
\end{equation}
where the computational domain is discretized into $M$ cells and $K$ is the grid cell number.
In the following sections, we compare four types of solution methods for the discrete two-phase problem \eqref{eq:definition_full_phase_based_residual}:
\begin{itemize}
\item FIM based on Newton's method with damping, in which the damping parameter is chosen to ensure that the largest saturation change between two Newton iterations is smaller than 0.2. The residual is computed by discretizing \eqref{eq:TwoPhaseFlowStandardForm} and using PPU to approximate the mobilities.
\item The sequential fully implicit method \cite{aMassConservativSequentialImplicitMultiscaleMethodForIsothermalEquationOfStateCompositionalProblems_MoynerTchelepi,aMultiscaleRestrictionSmoothedBasisMethodForCompositionalModels_MoynerTchelepi,ConsistentUpwindingForSFImultiscaleCompositionalSimulation_MoncorgeMoynerTchelepi,SuccessfulApplicationOfMultiscaleMethodsInArealReservoirSimulatorEnvironment_LieMoynerNatvig_et_al,SFIformulationForCompositionalSimulationUsingNaturalVariables_MoncorgeTchelepiJenny} referred to as SFI-$u_T$. The outer iteration consists in two steps. We first solve the flow problem nonlinearly and compute a new total velocity field using the updated pressure. In a second step, we solve the transport problem nonlinearly with a fixed total velocity. The residual is computed by discretizing \eqref{eq:fracFlowFormulation} and using IHU to approximate the mobilities. We have observed very slow convergence rates for SFI-$p$---in which the pressure is fixed during the transport solve---and we therefore do not report these results in the next sections.
\item The FSMSN method of Algorithm~\ref{algo:FSMSN}. We consider two versions of the algorithm. In FSMSN-$p$, the transport problem is computed with a fixed pressure. The residual is computed by discretizing \eqref{eq:TwoPhaseFlowStandardForm} and using PPU to approximate the mobilities. In FSMSN-$u_T$, the transport problem is computed with a fixed total velocity. The residual is computed by discretizing \eqref{eq:fracFlowFormulation} and using IHU to approximate the mobilities.
\item The MSPIN method of Algorithm~\ref{algo:mspin} referred to as MSPIN-$p$. The transport problem is computed with a fixed pressure. The residual is computed by discretizing \eqref{eq:TwoPhaseFlowStandardForm} and using PPU to approximate the mobilities.
\end{itemize}

The algorithms have been implemented in Matlab (R2019b) and tested on a basic laptop (Intel Core i5-8250U QuadCore @ 1.60GHz; 8GB RAM; hard disk size 256GB; Windows 10). The linear solves required by the algorithms have been performed by the default direct solver in Matlab (the "backslash" operator). No external library has been used.

\subsection{SPE10 bottom layer}
\label{subsec:spe10_bottom_layer}

\begin{figure}[!h]
\centering
\subfigure[Permeability (log(mD))]{
\includegraphics[height=0.3\textwidth]{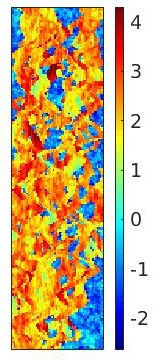}
\label{fig:permeabilityMap_spe10_bottom_layer.pdf}
}
\hspace{1.0cm}
\subfigure[Water saturation (0.01 PVI)]{
\includegraphics[height=0.3\textwidth]{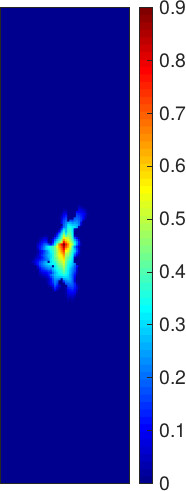}
\label{fig:spe10_bottom_layer_wettingPhase_t50}
}
\hspace{1.0cm}
\subfigure[Water saturation (0.05 PVI)]{
\includegraphics[height=0.3\textwidth]{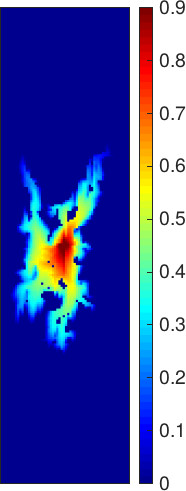}
\label{fig:spe10_bottom_layer_wettingPhase_t250}
}
\hspace{1.0cm}
\subfigure[Water saturation (0.1 PVI)]{
\includegraphics[height=0.3\textwidth]{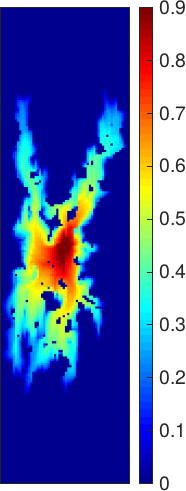}
\label{fig:spe10_bottom_layer_wettingPhase_t500}
}
\caption{SPE10 bottom layer: permeability map and water saturation maps at various times.}
\label{fig:spe10_bottom_layer_permeability_map_and_wetting_phase}
\end{figure}

We first consider a horizontal test case consisting of 13,200 cells in which the porosity and permeability fields are taken from the bottom layer of the SPE10 test case \cite{spe10}.
We inject the wetting phase (water) from the middle well and produce from the wells located in the four corners.
Capillary pressure is neglected.
The phase densities are set to $\rho_w = 1 025.0$ kg.m$^{-3}$ and $\rho_{nw} = 849.0$ kg.m$^{-3}$ and
the phase viscosities are set to $\mu_w=0.0003$ Pa.s and $\mu_{nw} = 0.003$ Pa.s.
We use quadratic Corey-type relative permeabilities.
The domain is initially fully saturated with the non-wetting phase.
We simulate 500 days of injection (0.1 total pore volume injected) with a constant time step size.
The permeability and final saturation maps are shown in Fig.~\ref{fig:spe10_bottom_layer_permeability_map_and_wetting_phase}. %\textcolor{red}{TODO: Use the same colorbar in Figs.1b-1c-1d}.

The nonlinear behavior of the schemes is illustrated in Fig.~\ref{fig:cv_curve_spe10_bottom_layer} for two selected time steps and in Table~\ref{tab:SPE10_bottom_layer_cfl69} for the full simulation.
We consider first a fixed tolerance of $10^{-6}$ in the pressure and transport subproblems.
We observe that for FSMSN-$u_T$, FSMSN-$p$, and MSPIN-$p$, the residual norm decreases at a faster rate than with Newton-based FIM and with SFI.
This faster rate results in a significant reduction in the number of outer iterations.
We note, however, that this reduction in the number of outer iterations requires performing a large number of subproblem iterations in FSMSN and MSPIN, which is likely to make these algorithms unpractical.
For instance, FSMSN-$u_T$ requires only 48 outer iterations---while 149 iterations are performed by Newton's method---but involves 81 pressure iterations and 205 transport iterations.

\pgfplotstableread
{
% cv curves obtained with normalized L2-norm of the residual
colnames     0  1  2  3  4  5  6  7  8  9  10
iter         1  2  3  4  5  6  7  8  9  10  11
fsmsn_ut     3.297010e-01  3.547662e-04  2.474423e-10  nan  nan  nan  nan  nan  nan  nan  nan
fsmsn_p      2.407295e+00  3.046874e-01  8.475359e-03  7.723031e-06  2.267107e-10  nan  nan  nan  nan  nan nan
mspin_p      2.693594e+00  3.781668e-01  1.152905e-02  1.538782e-05  2.729957e-10  nan  nan  nan  nan  nan nan
newton       4.483828e+00  3.376141e+00  1.130886e+00  1.934559e-01  2.940110e-02  2.095540e-03  3.330946e-05  1.488144e-08  nan  nan nan
sfi_ut       2.239684e+01  7.890311e-01  4.027486e-02  3.777096e-03  6.474023e-04  1.218995e-04  2.362014e-05  4.504522e-06  8.568434e-07  1.623822e-07  3.074351e-08
% cv curves obtained with L-infty morm of the residual
%colnames     0  1  2  3  4  5  6  7  8  9  10  11  12  13
%iter	      1  2  3  4  5  6  7  8  9  10  11  12  13  14
%fsmsn_ut     3.920600015  0.006324370067  2.082933381e-09  nan nan  nan  nan  nan  nan  nan   nan  nan  nan  nan
%fsmsn_p      39.389957  6.060707748  0.1758995594  0.000147875533  1.497615532e-09    nan  nan  nan  nan  nan   nan  nan  nan  nan
%mspin_p      38.91733534  6.770363794  0.2438218952  0.0003114332746  1.743817438e-09  nan  nan  nan  nan  nan   nan  nan  nan  nan
%newton       39.08243864  33.11047674  9.69875934  1.77924927  0.3533751226  0.02994786655  0.00040897783  2.103747681e-07  1.497615532e-09  nan  nan  nan  nan  nan
%sfi_ut       461.8744918  13.28388865  0.4634595912  0.04683253385  0.01019020678  0.002007359011  0.0003976757764  7.619561417e-05  1.455095276e-05  2.761754606e-06  5.234505069e-07  9.907637599e-08  1.87465119e-08  3.541584359e-09
}\CVCurveSecondTimestepSPETenBottomLayerData
\pgfplotstabletranspose[colnames from=colnames]\CVCurveSecondTimestepSPETenBottomLayerData{\CVCurveSecondTimestepSPETenBottomLayerData}

\pgfplotstableread
{
% cv curves obtained with normalized L2-norm of the residual
colnames     0  1  2  3  4  5  6  7
iter         1  2  3  4  5  6  7  8
fsmsn_ut     6.633197e-02  1.204729e-04  2.186475e-10  nan  nan  nan  nan  nan
fsmsn_p      6.768360e-02  1.540243e-04  2.354416e-10  nan  nan  nan  nan  nan
mspin_p      1.898311e-01  3.518014e-03  1.199781e-06  3.105109e-10  nan  nan  nan  nan
newton       6.290537e-01  2.297632e-01  2.536657e-02  5.757048e-04  7.092766e-07  2.256964e-10  nan  nan
sfi_ut       7.376366e-01  7.850188e-02  5.454349e-03  4.441061e-04  4.014570e-05  4.463489e-06  5.074512e-07  5.560213e-08
% cv curves obtained with L-infty morm of the residual
%colnames     0  1  2  3  4  5  6  7  8  9  10
%iter	      1  2  3  4  5  6  7  8  9  10  11
%fsmsn_ut     0.8059861008  0.002158776483  1.109136027e-09   nan  nan  nan  nan  nan  nan  nan  nan
%fsmsn_p      0.8071804466  0.002648009338  1.473349998e-09   nan  nan  nan  nan  nan  nan  nan  nan
%mspin_p      0.9498924245  0.05196143165  1.456674682e-05  1.033991059e-09 nan  nan  nan  nan  nan  nan  nan
%newton       4.433998829  3.092476941  0.2939146115  0.004594633988  6.411634506e-06  1.417741126e-09  nan  nan  nan  nan  nan
%sfi_ut       2.097765989  0.7841301083  0.06187748579  0.0052726231  0.0002188937814  4.232378991e-05  4.988109543e-06  4.581987396e-07  4.420050131e-08  1.230705493e-08  9.038860434e-09
}\CVCurveLastTimestepSPETenBottomLayerData
\pgfplotstabletranspose[colnames from=colnames]\CVCurveLastTimestepSPETenBottomLayerData{\CVCurveLastTimestepSPETenBottomLayerData}

\begin{figure}[h!]
  \centering
  \subfigure[Time step from 25 to 50 days]{
  \begin{tikzpicture}
    \begin{semilogyaxis}[%loglogaxis}[
      width = 0.4\textwidth,
      height = 0.4\textwidth,
      grid = major,
      major grid style= {very thin,draw=gray!30},
      xmin=0,xmax=12,
      ymin=1e-11,ymax=100,
      xlabel={Outer iterations},
      ylabel={Normalized $\ell_2$-norm of the residual},
      xtick={1,2,3,4,5,6,7,8,9,10,11},
      xticklabels={1,2,3,4,5,6,7,8,9,10,11},      
      ytick={1e-8, 1e-4, 1},
      yticklabels={$10^{-8}$,$10^{-4}$,$ 1$},
      %ytick distance=10^2,
      ylabel near ticks,
      xlabel near ticks,
      legend style={font=\small},
      tick label style={font=\small},
      %xticklabels from table={\CFLSPETenTopLayerHighInjectionData}{cfl_labels},
      log ticks with fixed point,
      label style={font=\small},
      %legend entries={ Newton, SFI-$u_T$, FSMSN-$u_T$, MSPIN-$p$, FSMSN-$p$},
      %legend cell align={left},
      %legend pos=outer north east,
      ]
      \addplot [mark=x, black, very thick] table [x=iter, y=newton] {\CVCurveSecondTimestepSPETenBottomLayerData};
      \addplot [mark=diamond, scarletred1, very thick] table [x=iter, y=sfi_ut] {\CVCurveSecondTimestepSPETenBottomLayerData};      
      \addplot [mark=square, skyblue1, very thick] table [x=iter, y=fsmsn_ut] {\CVCurveSecondTimestepSPETenBottomLayerData};
      \addplot [mark=asterisk, gray, very thick] table [x=iter, y=mspin_p] {\CVCurveSecondTimestepSPETenBottomLayerData};      
      \addplot [mark=o, chameleon1, very thick] table [x=iter, y=fsmsn_p] {\CVCurveSecondTimestepSPETenBottomLayerData};
    \end{semilogyaxis}
  \end{tikzpicture}
  \label{fig:cv_curve_secondTimestep_spe10_BottomLayer}
  }\hspace{0.25cm}
  \subfigure[Time step from 475 to 500 days]{
    \begin{tikzpicture}
    \begin{semilogyaxis}[
      width = 0.4\textwidth,
      height = 0.4\textwidth,
      grid = major,
      major grid style= {very thin,draw=gray!30},
      xmin=0,xmax=9,
      ymin=1e-11,ymax=100,
      xlabel={Outer iterations},
      ylabel={},
      xtick={1,2,3,4,5,6,7,8},
      xticklabels={1,2,3,4,5,6,7,8},      
      ytick={1e-8, 1e-4,1},
      yticklabels={$10^{-8}$,$10^{-4}$,$1$},
      %ytick distance=10^2,
      ylabel near ticks,
      xlabel near ticks,
      legend style={font=\small},
      tick label style={font=\small},
      %xticklabels from table={\CVCurveLastTimestepSPETenBottomLayerData}{cfl_labels},
      log ticks with fixed point,
      label style={font=\small},
      legend entries={ Newton, SFI-$u_T$, FSMSN-$u_T$, MSPIN-$p$, FSMSN-$p$},
      legend cell align={left},
      legend pos=outer north east,
      ]
      \addplot [mark=x, black, very thick] table [x=iter, y=newton] {\CVCurveLastTimestepSPETenBottomLayerData};
      \addplot [mark=diamond, scarletred1, very thick] table [x=iter, y=sfi_ut] {\CVCurveLastTimestepSPETenBottomLayerData};      
      \addplot [mark=square, skyblue1, very thick] table [x=iter, y=fsmsn_ut] {\CVCurveLastTimestepSPETenBottomLayerData};
      \addplot [mark=asterisk, gray, very thick] table [x=iter, y=mspin_p] {\CVCurveLastTimestepSPETenBottomLayerData};      
      \addplot [mark=o, chameleon1, very thick] table [x=iter, y=fsmsn_p] {\CVCurveLastTimestepSPETenBottomLayerData};
    \end{semilogyaxis}
    \end{tikzpicture}
    \label{fig:fig:cv_curve_lastTimestep_spe10_BottomLayer}
    }
  \caption{SPE10 bottom layer: Residual norm as a function of the number of outer iterations for two time steps. The plot on the left (respectively, on the right) corresponds to the second time step (respectively, the last time step) in the simulation. The max CFL number for these two time steps is approximately 69.}
  \label{fig:cv_curve_spe10_bottom_layer}
\end{figure}
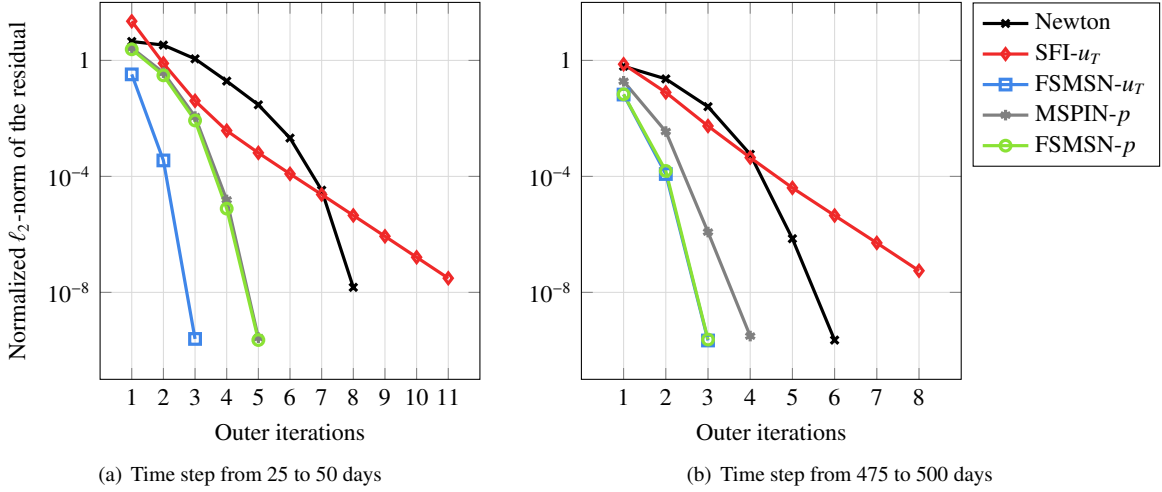
\begin{table}[!h]
  \centering
  \small
  \begingroup
  \setlength{\tabcolsep}{6pt}
  \begin{tabular}{l     c      c      c     c     c     c     c}
    \toprule
    Solver &                            && \multicolumn{2}{c}{Fixed total velocity}                        && \multicolumn{2}{c}{Fixed pressure} \\
    \cmidrule{4-5} \cmidrule{7-8} 
           & \multicolumn{1}{c}{Newton} && \multicolumn{1}{c}{SFI-$u_T$} & \multicolumn{1}{c}{FSMSN-$u_T$} && \multicolumn{1}{c}{MSPIN-$p$} & \multicolumn{1}{c}{FSMSN-$p$} \\
    \midrule
    Nonlinear iterations      & 149                   && 132 & 48  && 68 & 56 \\
    Iterations per time step  & 7.5                   && 6.6 & 2.4 && 3.4 & 2.8 \\
    Pressure iterations       & \multicolumn{1}{c}{-} && 204 & 81  && 116 & 93 \\
    Transport iterations      & \multicolumn{1}{c}{-} && 328 & 205 && 248 & 215 \\
    \bottomrule
  \end{tabular}
  \endgroup
\caption{SPE10 bottom layer: nonlinear behavior of the schemes with a fixed subproblem tolerance of $10^{-6}$. The maximum CFL number is 69 and the total PVI is 0.1.}
\label{tab:SPE10_bottom_layer_cfl69}
\end{table}
\begin{table}[h!]
  \centering
  \small
  \begingroup
  \setlength{\tabcolsep}{6pt}
  \begin{tabular}{l      j      j      j      c      j      j      j      c      j      j      j      c      j      j      j} 
    \toprule
    Solver                    & \multicolumn{7}{c}{Fixed total velocity}                         && \multicolumn{7}{c}{Fixed pressure} \\
    \cmidrule{2-8} \cmidrule{10-16} 
                              & \multicolumn{3}{c}{SFI-$u_T$} && \multicolumn{3}{c}{FSMSN-$u_T$} && \multicolumn{3}{c}{MSPIN-$p$} && \multicolumn{3}{c}{FSMSN-$p$} \\
    \cmidrule{2-4}  \cmidrule{6-8}  \cmidrule{10-12}  \cmidrule{14-16}
                              & $A_1$ & $A_2$ & $A_3$   && $A_1$ & $A_2$ & $A_3$     && $A_1$ & $A_2$ & $A_3$   && $A_1$ & $A_2$ & $A_3$   \\
    \midrule
    Nonlinear iterations      & 133 & 136 & 133 && 62  & 66  & 63  && 96  & 114 & 102 && 66  & 71 & 68 \\
    Iterations per time step  & 6.7 & 6.8 & 6.7 && 3.1 & 3.3 & 3.2 && 4.8 & 5.7 & 5.1 && 3.3 & 3.6 & 3.4 \\
    Pressure iterations       & 134 & 137 & 135 && 63  & 67  & 64  && 97  & 115 & 103 && 67  & 72 & 69 \\
    Transport iterations      & 188 & 178 & 180 && 109 & 94  & 107 && 152 & 150 & 151 && 110 & 99 & 110\\
    \bottomrule
  \end{tabular}
  \endgroup
  \caption{SPE10 bottom layer: nonlinear behavior of the schemes with an adaptive subproblem tolerance computed using the strategies $(A_1)$, $(A_2)$, and $(A_3)$ of Section~\ref{sec:adaptive_tolerance}. The maximum CFL number is 69 and the total PVI is 0.1.} 
\label{tab:adaptiveTolerance_SPE10_bottom_layer_cfl69}
\end{table}

Table~\ref{tab:adaptiveTolerance_SPE10_bottom_layer_cfl69} shows the results of the adaptive strategies introduced in Section~\ref{sec:adaptive_tolerance} to relax the tolerance in the subproblems and reduce the computational cost of the pressure and transport steps in FSMSN and MSPIN. 
For both FSMSN and MSPIN, using adaptive tolerances in the subproblems results in a slight increase in the number of outer iterations, while the number of pressure iterations and---more significantly---the number of transport iterations are reduced.

\begin{figure}[h!]
  \centering
  \subfigure[Fixed tolerance]{
  \begin{tikzpicture}
    \begin{loglogaxis}[
      width = 0.4\textwidth,
      height = 0.4\textwidth,
      grid = major,
      major grid style= {very thin,draw=gray!30},
      xmin=20,xmax=380,
      ymin=10,ymax=400,
      xlabel={Maximum CFL number},
      ylabel={Cumulative number of outer iterations},
      xtick=data, 
      ytick={10,20,40,100,200,400},
      yticklabels={10,20,40,100,200,400},
      ytick distance=10^1,
      ylabel near ticks,
      xlabel near ticks,
      legend style={font=\small},
      tick label style={font=\small},
      %xticklabels from table={\CFLSPETenBottomLayerData}{cfl_labels},
      log ticks with fixed point,
      label style={font=\small},
      %legend entries={ Newton, SFI-$u_T$, FSMSN-$u_T$, MSPIN-$p$, FSMSN-$p$},
      %legend cell align={left},
      %legend pos=outer north east,
      ]
      \addplot [mark=x, black, very thick] table [x=cfl, y=newton] {\CFLSPETenBottomLayerData};
      \addplot [mark=diamond, scarletred1, very thick] table [x=cfl, y=sfi_ut] {\CFLSPETenBottomLayerData};      
      \addplot [mark=square, skyblue1, very thick] table [x=cfl, y=fsmsn_ut] {\CFLSPETenBottomLayerData};
      \addplot [mark=asterisk, gray, very thick] table [x=cfl, y=mspin_p] {\CFLSPETenBottomLayerData};      
      \addplot [mark=o, chameleon1, very thick] table [x=cfl, y=fsmsn_p] {\CFLSPETenBottomLayerData};
    \end{loglogaxis}
  \end{tikzpicture}
  \label{fig:spe10_bottom_layer_cumulative_iter_vs_cfl_fixed_tolerance}
  }\hspace{0.5cm}
  \subfigure[Adaptive tolerance with $(A_1)$]{
  \begin{tikzpicture}
    \begin{loglogaxis}[
      width = 0.4\textwidth,
      height = 0.4\textwidth,
      grid = major,
      major grid style= {very thin,draw=gray!30},
      xmin=20,xmax=380,
      ymin=10,ymax=400,
      xlabel={Maximum CFL number},
      ylabel={Cumulative number of outer iterations},
      xtick=data, 
      ytick={10,20,40,100,200,400},
      yticklabels={10,20,40,100,200,400},
      ytick distance=10^1,
      ylabel near ticks,
      xlabel near ticks,
      legend style={font=\small},
      tick label style={font=\small},
      %xticklabels from table={\CFLSPETenBottomLayerAdaptivetoleranceData}{cfl_labels},
      log ticks with fixed point,
      label style={font=\small},
      legend entries={ Newton, SFI-$u_T$, FSMSN-$u_T$, MSPIN-$p$, FSMSN-$p$},
      legend cell align={left},
      legend pos=outer north east,
      ]
      \addplot [mark=x, black, very thick] table [x=cfl, y=newton] {\CFLSPETenBottomLayerAdaptivetoleranceData};
      \addplot [mark=diamond, scarletred1, very thick] table [x=cfl, y=sfi_ut] {\CFLSPETenBottomLayerAdaptivetoleranceData};
      \addplot [mark=square, skyblue1, very thick] table [x=cfl, y=fsmsn_ut] {\CFLSPETenBottomLayerAdaptivetoleranceData};
      \addplot [mark=asterisk, gray, very thick] table [x=cfl, y=mspin_p] {\CFLSPETenBottomLayerAdaptivetoleranceData};
      \addplot [mark=o, chameleon1, very thick] table [x=cfl, y=fsmsn_p] {\CFLSPETenBottomLayerAdaptivetoleranceData};
    \end{loglogaxis}
  \end{tikzpicture}
  \label{fig:spe10_bottom_layer_cumulative_iter_vs_cfl_adaptive_tolerance}
  }
  \caption{SPE10 bottom layer: cumulative number of outer iterations for the full simulation as a function of the maximum CFL number observed during the simulation.}
  \label{fig:spe10_bottom_layer_cumulative_iter_vs_cfl}
\end{figure}
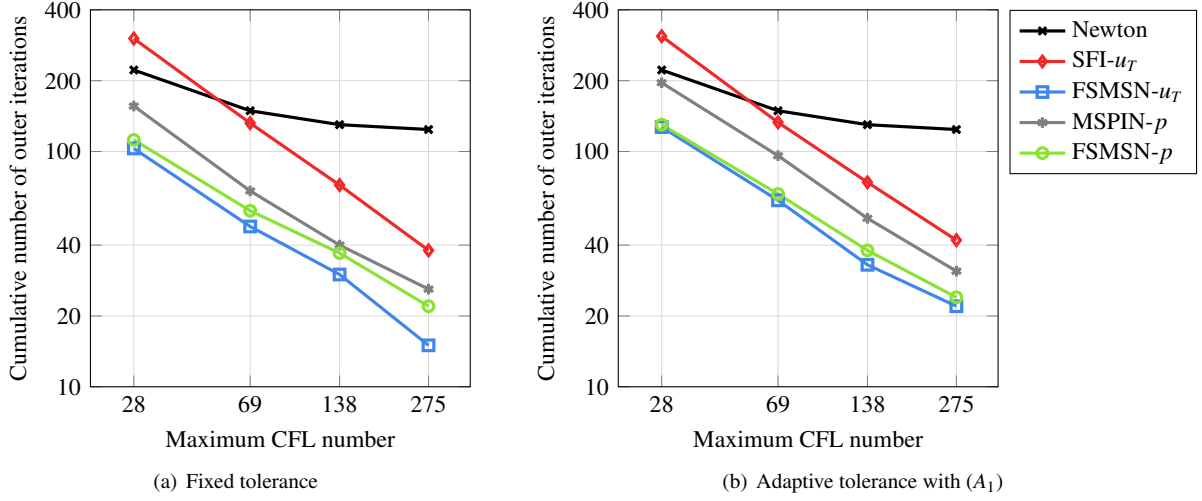

To conclude this section, we study the sensitivity of the results to the time step size---measured by the maximum CFL number observed during the simulation.
The results with a fixed tolerance and an adaptive tolerance are shown in Fig.~\ref{fig:spe10_bottom_layer_cumulative_iter_vs_cfl}.
We observe that when the time step is increased, the number of Newton iterations is slightly reduced for CFL numbers smaller than 138, but stagnates for CFL numbers larger than 138.
This is not the case for FSMSN-$u_T$, FSMSN-$p$, and MSPIN-$p$, as these schemes exhibit a significant reduction in the number of outer iterations as the time step is increased, even for large CFL numbers.
With both fixed tolerance and adaptive tolerance in the pressure and transport subproblems, FSMSN-$u_T$ is the solution strategy that requires the smallest number of outer iterations.

%In Fig.~\ref{fig:spe10_bottom_layer_cumulative_iter_vs_cfl}, we show the total number of nonlinear iterations with respect to the maximum CFL number for different timestep.  
%During one simulation, the timestep remains fixed for all algorithms.
\pgfplotstableread
{
colnames     0	   1      2      3
cfl	     28    69	  138    275
%cfl_labels   $10^{1.5}$   $10^2$   $10^{2.5}$   $10^{3}$
fsmsn_ut     103   48     30     15
fsmsn_p      112   56     37     22
mspin_p      156   68     40     26
newton       222   149    130    124
sfi_ut       302   132    72     38
}\CFLSPETenBottomLayerData
\pgfplotstabletranspose[colnames from=colnames]\CFLSPETenBottomLayerData{\CFLSPETenBottomLayerData}

\pgfplotstableread
{
colnames     0	   1      2      3
cfl	     28    69	  138    275
%cfl_labels   $10^{1.5}$   $10^2$   $10^{2.5}$   $10^{3}$
fsmsn_ut     127   62     33     22
fsmsn_p      130   66     38     24
mspin_p      196   96     52     31
newton       222   149    130    124
sfi_ut       309   133    74     42
}\CFLSPETenBottomLayerAdaptivetoleranceData
\pgfplotstabletranspose[colnames from=colnames]\CFLSPETenBottomLayerAdaptivetoleranceData{\CFLSPETenBottomLayerAdaptivetoleranceData}

\subsection{SPE10 top layer}
\label{subsec:spe10_top_layer}

In this example, we study the impact of the introduction of buoyancy forces on the nonlinear behavior of the schemes.
To do that, we consider a tilted two-dimensional domain in which the porosity and permeability fields are taken from the top layer of the SPE10 test case.
The fluid properties and the well locations are the same as in the previous section.
We inject 0.08 total pore volume in two distinct configurations:
\begin{itemize}
\item Case 1: tilting of $60^{\circ}$ in the $y$-direction and fast injection rate (9.352 m$^3$/day) %(58.82 bbl/day) 
so that viscous forces dominate.
\item Case 2: same tilting and slower injection rate (0.9352 m$^3$/day) %(5.882 bbl/day) 
so that buoyancy forces dominate.
\end{itemize}
We note that the time step size and total simulation time are adapted to obtain approximately the same maximum CFL number and total PVI in the two cases.
We start the simulation with a uniform initial saturation field equal to $S^{0}_w = 0$.

The nonlinear behavior of the schemes with a maximum CFL number of approximately 69 is summarized in Tables~\ref{tab:SPE10_top_layer_HighInjectionRate_cfl69}-\ref{tab:adaptiveTolerance_SPE10_topLayer_cfl69_highInjection} for Case 1, and in Tables~\ref{tab:SPE10_top_layer_lowInjectionRate_cfl69}-\ref{tab:adaptiveTolerance_SPE10_topLayer_cfl69_LowInjection} for Case 2.
We observe that increasing the strength of buoyancy forces (relatively to viscous forces) makes the nonlinear convergence of SFI-$u_T$ very slow, and deteriorates slightly that of MSPIN-$p$.
However, for Newton's method, FSMSM-$u_T$, and FSMSN-$p$, we observe a slight improvement in the nonlinear behavior when buoyancy forces are stronger.
\begin{table}[h!]
  \centering
  \small
  \begingroup
  \setlength{\tabcolsep}{6pt}
  \begin{tabular}{l     c      c      c     c     c     c     c}
    \toprule
    Solver &                            && \multicolumn{2}{c}{Fixed total velocity}                        && \multicolumn{2}{c}{Fixed pressure} \\
    \cmidrule{4-5} \cmidrule{7-8} 
           & \multicolumn{1}{c}{Newton} && \multicolumn{1}{c}{SFI-$u_T$} & \multicolumn{1}{c}{FSMSN-$u_T$} && \multicolumn{1}{c}{MSPIN-$p$} & \multicolumn{1}{c}{FSMSN-$p$} \\
    \midrule
    Nonlinear iterations      & 135                   && 130  & 51   && 60   & 56 \\
    Iterations per time step  & 8.43                  && 8.12 & 3.19 && 3.75 & 3.5 \\
    Pressure iterations       & \multicolumn{1}{c}{-} && 496  & 260  && 314  & 275 \\
    Transport iterations      & \multicolumn{1}{c}{-} && 288  & 173  && 201  & 180 \\
    \bottomrule
  \end{tabular}
  \endgroup
  \caption{SPE10 top layer: nonlinear behavior of the schemes with a fixed subproblem tolerance of $10^{-6}$ (Case 1 with fast injection rate). Approximately 2\% of the interfaces experience counter-current flow. The maximum CFL number is 69 and the total PVI is 0.08.}
\label{tab:SPE10_top_layer_HighInjectionRate_cfl69}
\end{table}
\begin{table}[h!]
  \centering
  \small
  \begingroup
  \setlength{\tabcolsep}{6pt}
  \begin{tabular}{l      j      j      j      c      j      j      j      c      j      j      j      c      j      j      j} 
    \toprule
    Solver                    & \multicolumn{7}{c}{Fixed total velocity}                         && \multicolumn{7}{c}{Fixed pressure} \\
    \cmidrule{2-8} \cmidrule{10-16} 
                              & \multicolumn{3}{c}{SFI-$u_T$} && \multicolumn{3}{c}{FSMSN-$u_T$} && \multicolumn{3}{c}{MSPIN-$p$} && \multicolumn{3}{c}{FSMSN-$p$} \\
    \cmidrule{2-4}  \cmidrule{6-8}  \cmidrule{10-12}  \cmidrule{14-16}
                              & $A_1$ & $A_2$ & $A_3$   && $A_1$ & $A_2$ & $A_3$     && $A_1$ & $A_2$ & $A_3$   && $A_1$ & $A_2$ & $A_3$ \\
    \midrule
    Nonlinear iterations      & 158 & d.n.c. & 164  && 51  & 56  & 55  && 80  & 88  & 89  && 56  & 58  & 56  \\
    Iterations per time step  & 9.9 & d.n.c. & 10.3 && 3.2 & 3.5 & 3.4 && 5   & 5.5 & 5.7 && 3.5 & 3.6 & 3.5 \\
    Pressure iterations       & 178 & d.n.c. & 174  && 69  & 79  & 69  && 109 & 114 & 110 && 74  & 77  & 73  \\
    Transport iterations      & 196 & d.n.c. & 198  && 85  & 80  & 90  && 122 & 115 & 131 && 90  & 83  & 96  \\
    \bottomrule
  \end{tabular}
  \endgroup
  \caption{SPE10 top layer: nonlinear behavior of the schemes with an adaptive subproblem tolerance computed using the strategies $(A_1)$, $(A_2)$, and $(A_3)$ of Section~\ref{sec:adaptive_tolerance} (Case 1 with fast injection rate). The maximum CFL number is 69 and the total PVI is 0.08. d.n.c. denotes lack of convergence.}
\label{tab:adaptiveTolerance_SPE10_topLayer_cfl69_highInjection}
\end{table}

\begin{table}[h!]
  \centering
  \small
  \begingroup
  \setlength{\tabcolsep}{6pt}
  \begin{tabular}{l     c      c      c     c     c     c     c}
    \toprule
    Solver &                            && \multicolumn{2}{c}{Fixed total velocity}                        && \multicolumn{2}{c}{Fixed pressure} \\
    \cmidrule{4-5} \cmidrule{7-8} 
           & \multicolumn{1}{c}{Newton} && \multicolumn{1}{c}{SFI-$u_T$} & \multicolumn{1}{c}{FSMSN-$u_T$} && \multicolumn{1}{c}{MSPIN-$p$} & \multicolumn{1}{c}{FSMSN-$p$} \\
    \midrule
    Nonlinear iterations      & 127                   && 567   & 48  && 67   & 56 \\
    Iterations per time step  & 7.93                  && 33.35 & 3   && 3.18 & 3.5 \\
    Pressure iterations       & \multicolumn{1}{c}{-} && 1049  & 86  && 120  & 96 \\
    Transport iterations      & \multicolumn{1}{c}{-} && 2166  & 188 && 244  & 214 \\
    \bottomrule
  \end{tabular}
  \endgroup
  \caption{SPE10 top layer: nonlinear behavior of the schemes with a fixed subproblem tolerance of $10^{-6}$ (Case 2 with a slower injection rate). Approximately 8\% of the interfaces experience counter-current flow. The maximum CFL number is 69 and the total PVI is 0.08.}
  \label{tab:SPE10_top_layer_lowInjectionRate_cfl69}
\end{table}
\begin{table}[h!]
  \centering
  \small
  \begingroup
  \setlength{\tabcolsep}{6pt}
  \begin{tabular}{l      j      j      j      c      j      j      j      c      j      j      j      c      j      j      j} 
    \toprule
    Solver                    & \multicolumn{7}{c}{Fixed total velocity}                         && \multicolumn{7}{c}{Fixed pressure} \\
    \cmidrule{2-8} \cmidrule{10-16} 
                              & \multicolumn{3}{c}{SFI-$u_T$} && \multicolumn{3}{c}{FSMSN-$u_T$} && \multicolumn{3}{c}{MSPIN-$p$} && \multicolumn{3}{c}{FSMSN-$p$} \\
    \cmidrule{2-4}  \cmidrule{6-8}  \cmidrule{10-12}  \cmidrule{14-16}
                              & $A_1$ & $A_2$ & $A_3$   && $A_1$ & $A_2$ & $A_3$     && $A_1$ & $A_2$ & $A_3$   && $A_1$ & $A_2$ & $A_3$ \\
    \midrule
    Nonlinear iterations      & 566  & d.n.c. & 571  && 49  & 49  & 51  && 79  & 93  & 91  && 55  & 55  & 56 \\
    Iterations per time step  & 33.3 & d.n.c. & 33.6 && 3.1 & 3.1 & 3.2 && 4.9 & 5.8 & 5.7 && 3.4 & 3.4 & 3.5 \\
    Pressure iterations       & 975  & d.n.c. & 981  && 51  & 51  & 54  && 84  & 97  & 94  && 58  & 58  & 60 \\
    Transport iterations      & 1229 & d.n.c. & 1141 && 89  & 93  & 90  && 130 & 129 & 129 && 92  & 95  & 92 \\
    \bottomrule
  \end{tabular}
  \endgroup
  \caption{SPE10 top layer: nonlinear behavior of the schemes with an adaptive subproblem tolerance computed using the strategies $(A_1)$, $(A_2)$, and $(A_3)$ of Section~\ref{sec:adaptive_tolerance} (Case 2 with slower injection rate). The maximum CFL number is 69 and the total PVI is 0.08. d.n.c. denotes lack of convergence.}
\label{tab:adaptiveTolerance_SPE10_topLayer_cfl69_LowInjection}
\end{table}

In this section, the nonlinear behavior of FSMSN-$u_T$ and FSMSN-$p$ with adaptive tolerance in the subproblems remains excellent.
For Cases 1 and 2, the reduction in the number of subproblem iterations obtained with the adaptive strategies is drastic but the increase in the number of outer iterations is small.
However, for MSPIN-$p$, we observe that using adaptive tolerances deteriorates quite significantly the nonlinear behavior, with for instance an increase by 48\% in the number of outer iterations for strategy $(A_3)$.

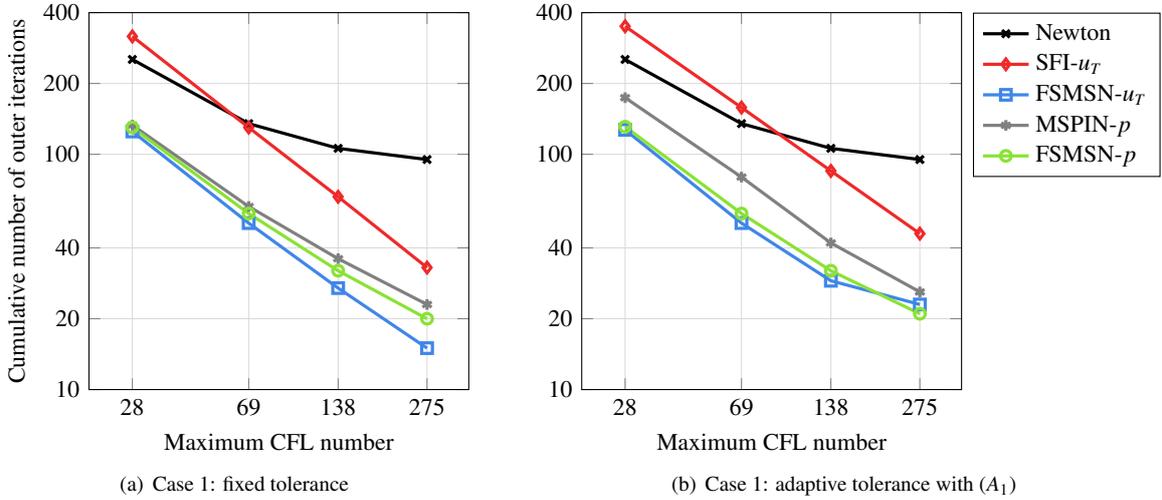
\begin{figure}[h!]
  \centering
  \subfigure[Case 1: fixed tolerance]{
  \begin{tikzpicture}
    \begin{loglogaxis}[
      width = 0.4\textwidth,
      height = 0.4\textwidth,
      grid = major,
      major grid style= {very thin,draw=gray!30},
      xmin=20,xmax=380,
      ymin=10,ymax=400,
      xlabel={Maximum CFL number},
      ylabel={Cumulative number of outer iterations},
      xtick=data,
      ytick={10,20,40,100,200,400},
      yticklabels={10,20,40,100,200,400},      
      ytick distance=10^1,
      ylabel near ticks,
      xlabel near ticks,
      legend style={font=\small},
      tick label style={font=\small},
      %xticklabels from table={\CFLSPETenTopLayerHighInjectionData}{cfl_labels},
      log ticks with fixed point,
      label style={font=\small},
      %legend entries={ Newton, SFI-$u_T$, FSMSN-$u_T$, MSPIN-$p$, FSMSN-$p$},
      %legend cell align={left},
      %legend pos=outer north east,
      ]
      \addplot [mark=x, black, very thick] table [x=cfl, y=newton] {\CFLSPETenTopLayerHighInjectionData};
      \addplot [mark=diamond, scarletred1, very thick] table [x=cfl, y=sfi_ut] {\CFLSPETenTopLayerHighInjectionData};      
      \addplot [mark=square, skyblue1, very thick] table [x=cfl, y=fsmsn_ut] {\CFLSPETenTopLayerHighInjectionData};
      \addplot [mark=asterisk, gray, very thick] table [x=cfl, y=mspin_p] {\CFLSPETenTopLayerHighInjectionData};      
      \addplot [mark=o, chameleon1, very thick] table [x=cfl, y=fsmsn_p] {\CFLSPETenTopLayerHighInjectionData};
    \end{loglogaxis}
  \end{tikzpicture}
  \label{fig:nonlinearIterCurve_spe10_top_layerHighInjection.pdf}
  }\hspace{0.5cm}
  \subfigure[Case 1: adaptive tolerance with $(A_1)$]{
    \begin{tikzpicture}
    \begin{loglogaxis}[
      width = 0.4\textwidth,
      height = 0.4\textwidth,
      grid = major,
      major grid style= {very thin,draw=gray!30},
      xmin=20,xmax=380,
      ymin=10,ymax=400,
      xlabel={Maximum CFL number},
      ylabel={},
      xtick=data,
      ytick={10,20,40,100,200,400},
      yticklabels={10,20,40,100,200,400},      
      ylabel near ticks,
      xlabel near ticks,
      legend style={font=\small},
      tick label style={font=\small},
      %xticklabels from table={\CFLSPETenTopLayerLowInjectionData}{cfl_labels},
      log ticks with fixed point,
      label style={font=\small},
      legend entries={ Newton, SFI-$u_T$, FSMSN-$u_T$, MSPIN-$p$, FSMSN-$p$},
      legend cell align={left},
      legend pos=outer north east,
      ]
      \addplot [mark=x, black, very thick] table [x=cfl, y=newton] {\CFLSPETenTopLayerHighInjectionAdaptivetoleranceData};
      \addplot [mark=diamond, scarletred1, very thick] table [x=cfl, y=sfi_ut] {\CFLSPETenTopLayerHighInjectionAdaptivetoleranceData};
      \addplot [mark=square, skyblue1, very thick] table [x=cfl, y=fsmsn_ut] {\CFLSPETenTopLayerHighInjectionAdaptivetoleranceData};
      \addplot [mark=asterisk, gray, very thick] table [x=cfl, y=mspin_p] {\CFLSPETenTopLayerHighInjectionAdaptivetoleranceData};
      \addplot [mark=o, chameleon1, very thick] table [x=cfl, y=fsmsn_p] {\CFLSPETenTopLayerHighInjectionAdaptivetoleranceData};
    \end{loglogaxis}
    \end{tikzpicture}
    \label{fig:nonlinearIterCurve_spe10_top_layerLowInjection.pdf}
  }
  \caption{SPE10 top layer: cumulative number of outer iterations for the full simulation as a function of the maximum CFL number observed during the simulation (Case 1 with strong viscous forces relatively to buoyancy forces).}
  \label{fig:nonlinearIterCurve_spe10_top_layer}
\end{figure}

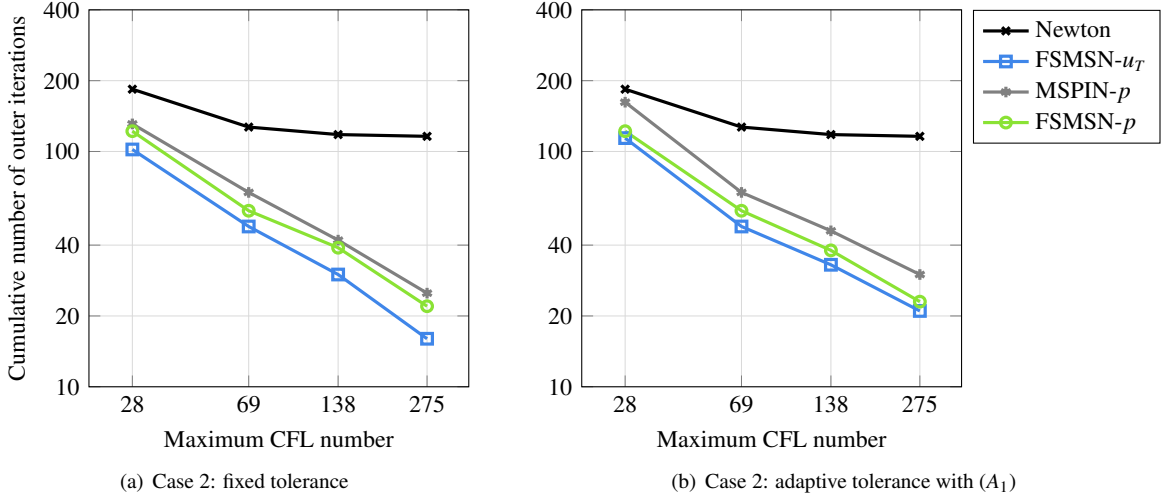
\begin{figure}[h!]
  \centering
  \subfigure[Case 2: fixed tolerance]{
  \begin{tikzpicture}
    \begin{loglogaxis}[
      width = 0.4\textwidth,
      height = 0.4\textwidth,
      grid = major,
      major grid style= {very thin,draw=gray!30},
      xmin=20,xmax=380,
      ymin=10,ymax=400,
      xlabel={Maximum CFL number},
      ylabel={Cumulative number of outer iterations},
      xtick=data,
      ytick={10,20,40,100,200,400},
      yticklabels={10,20,40,100,200,400},      
      ytick distance=10^1,
      ylabel near ticks,
      xlabel near ticks,
      legend style={font=\small},
      tick label style={font=\small},
      %xticklabels from table={\CFLSPETenTopLayerHighInjectionAdaptivetoleranceData}{cfl_labels},
      log ticks with fixed point,
      label style={font=\small},
      %legend entries={ Newton, SFI-$u_T$, FSMSN-$u_T$, MSPIN-$p$, FSMSN-$p$},
      %legend cell align={left},
      %legend pos=outer north east,
      ]
      \addplot [mark=x, black, very thick] table [x=cfl, y=newton] {\CFLSPETenTopLayerLowInjectionData};
      %\addplot [mark=diamond, scarletred1, very thick] table [x=cfl, y=sfi_ut] {\CFLSPETenTopLayerLowInjectionData};      
      \addplot [mark=square, skyblue1, very thick] table [x=cfl, y=fsmsn_ut] {\CFLSPETenTopLayerLowInjectionData};
      \addplot [mark=asterisk, gray, very thick] table [x=cfl, y=mspin_p] {\CFLSPETenTopLayerLowInjectionData};      
      \addplot [mark=o, chameleon1, very thick] table [x=cfl, y=fsmsn_p] {\CFLSPETenTopLayerLowInjectionData};      
    \end{loglogaxis}
  \end{tikzpicture}
  \label{fig:nonlinearIterCurve_spe10_top_layerHighInjection_adaptive_tolerance}
  }\hspace{0.5cm}
  \subfigure[Case 2: adaptive tolerance with $(A_1)$]{
    \begin{tikzpicture}
    \begin{loglogaxis}[
      width = 0.4\textwidth,
      height = 0.4\textwidth,
      grid = major,
      major grid style= {very thin,draw=gray!30},
      xmin=20,xmax=380,
      ymin=10,ymax=400,
      xlabel={Maximum CFL number},
      ylabel={},
      xtick=data,
      ytick={10,20,40,100,200,400},
      yticklabels={10,20,40,100,200,400},      
      ylabel near ticks,
      xlabel near ticks,
      legend style={font=\small},
      tick label style={font=\small},
      %xticklabels from table={\CFLSPETenTopLayerLowInjectionAdaptivetoleranceData}{cfl_labels},
      log ticks with fixed point,
      label style={font=\small},
      legend entries={ Newton, FSMSN-$u_T$, MSPIN-$p$, FSMSN-$p$},
      legend cell align={left},
      legend pos=outer north east,
      ]
      \addplot [mark=x, black, very thick] table [x=cfl, y=newton] {\CFLSPETenTopLayerLowInjectionAdaptivetoleranceData};
      %\addplot [mark=diamond, sacarletred1, very thick] table [x=cfl, y=sfi_ut] {\CFLSPETenTopLayerLowInjectionAdaptivetoleranceData};
      \addplot [mark=square, skyblue1, very thick] table [x=cfl, y=fsmsn_ut] {\CFLSPETenTopLayerLowInjectionAdaptivetoleranceData};
      \addplot [mark=asterisk, gray, very thick] table [x=cfl, y=mspin_p] {\CFLSPETenTopLayerLowInjectionAdaptivetoleranceData};
      \addplot [mark=o, chameleon1, very thick] table [x=cfl, y=fsmsn_p] {\CFLSPETenTopLayerLowInjectionAdaptivetoleranceData};
    \end{loglogaxis}
    \end{tikzpicture}
    \label{fig:nonlinearIterCurve_spe10_top_layerLowInjection_adaptive_tolerance}
  }
  \caption{SPE10 top layer: cumulative number of outer iterations for the full simulation as a function of the maximum CFL number observed during the simulation (Case 2 with strong buoyancy forces relatively to viscous forces).}
  \label{fig:nonlinearIterCurve_spe10_top_layer_with_adaptive_tolerance}
\end{figure}

The sensitivity of the nonlinear behavior to the time step size is studied in Figs.~\ref{fig:nonlinearIterCurve_spe10_top_layer} and \ref{fig:nonlinearIterCurve_spe10_top_layer_with_adaptive_tolerance}.
We note that these figures do not include the results obtained with SFI-$u_T$ as this approach does not converge for Case 2 when the CFL number is increased.
As in Section \ref{subsec:spe10_bottom_layer}, the number of Newton iterations stagnates for CFL numbers larger than 69, while the number of outer iterations required by FSMSN-$u_T$, FSMSN-$p$, and MSPIN-$p$ reduce for the range of time step sizes considered here.

\pgfplotstableread
{
colnames     0	   1      2      3
cfl	     28    69	  138    275
%cfl_labels   $10^{1.5}$   $10^2$   $10^{2.5}$   $10^{3}$
fsmsn_ut     125   51     27     15
fsmsn_p      130   56     32     20
mspin_p      133   60     36     23
newton       253   135    106    95
sfi_ut       317   130    66     33
}\CFLSPETenTopLayerHighInjectionData
\pgfplotstabletranspose[colnames from=colnames]\CFLSPETenTopLayerHighInjectionData{\CFLSPETenTopLayerHighInjectionData}

\pgfplotstableread
{
colnames     0	   1      2      3
cfl	     28    69	  138    275
%cfl_labels   $10^{1.5}$   $10^2$   $10^{2.5}$   $10^{3}$
fsmsn_ut     127   51     29     23
fsmsn_p      131   56     32     21
mspin_p      174   80     42     26
newton       253   135    106    95
sfi_ut       350   158    85     46
}\CFLSPETenTopLayerHighInjectionAdaptivetoleranceData
\pgfplotstabletranspose[colnames from=colnames]\CFLSPETenTopLayerHighInjectionAdaptivetoleranceData{\CFLSPETenTopLayerHighInjectionAdaptivetoleranceData}

\pgfplotstableread
{
colnames     0	   1      2      3
cfl	     28    69	  138    275
%cfl_labels   $10^{1.5}$   $10^2$   $10^{2.5}$   $10^{3}$
fsmsn_ut     102   48     30     16
fsmsn_p      122   56     39     22
mspin_p      131   67     42     25
newton       184   127    118    116
sfi_ut       505   567    nan    nan
}\CFLSPETenTopLayerLowInjectionData
\pgfplotstabletranspose[colnames from=colnames]\CFLSPETenTopLayerLowInjectionData{\CFLSPETenTopLayerLowInjectionData}

\pgfplotstableread
{
colnames     0	   1      2      3
cfl	     28    69	  138    275
%cfl_labels   $10^{1.5}$   $10^2$   $10^{2.5}$   $10^{3}$
fsmsn_ut     114   48     33     21
fsmsn_p      122   56     38     23
mspin_p      162   67     46     30
newton       184   127    118    116
sfi_ut       507   567    nan    nan
}\CFLSPETenTopLayerLowInjectionAdaptivetoleranceData
\pgfplotstabletranspose[colnames from=colnames]\CFLSPETenTopLayerLowInjectionAdaptivetoleranceData{\CFLSPETenTopLayerLowInjectionAdaptivetoleranceData}

\subsection{Gravity segregation test case}
\label{sec:lock_exchange}
\begin{figure}[!h]
\centering
\subfigure[Water saturation after 75 days]{
\includegraphics[height=0.27\textwidth]{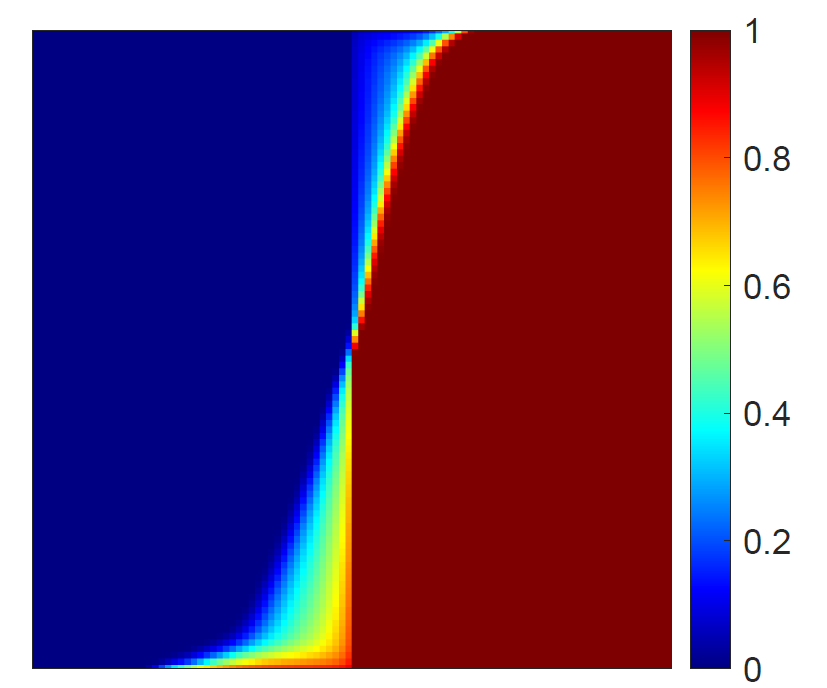}
\label{fig:lockExchangeHomoRelPerm_wetting_phase_t75}
}
\hspace{1.0cm}
\subfigure[Water saturation after 150 days]{
\includegraphics[height=0.27\textwidth]{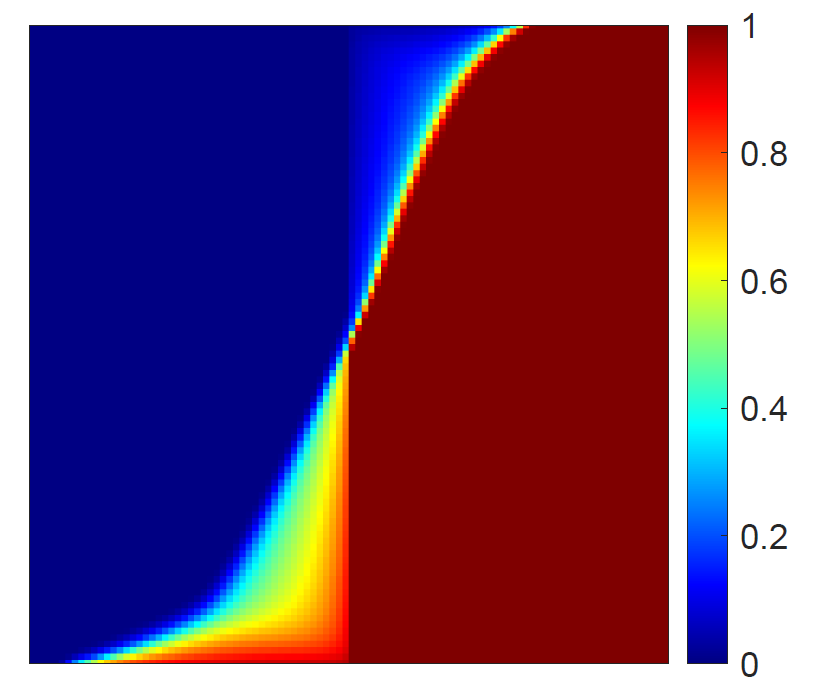}
\label{fig:lockExchangeHomoRelPerm_wetting_phase_t150}
}
\hspace{1.0cm}
\subfigure[Water saturation after 300 days]{
\includegraphics[height=0.27\textwidth]{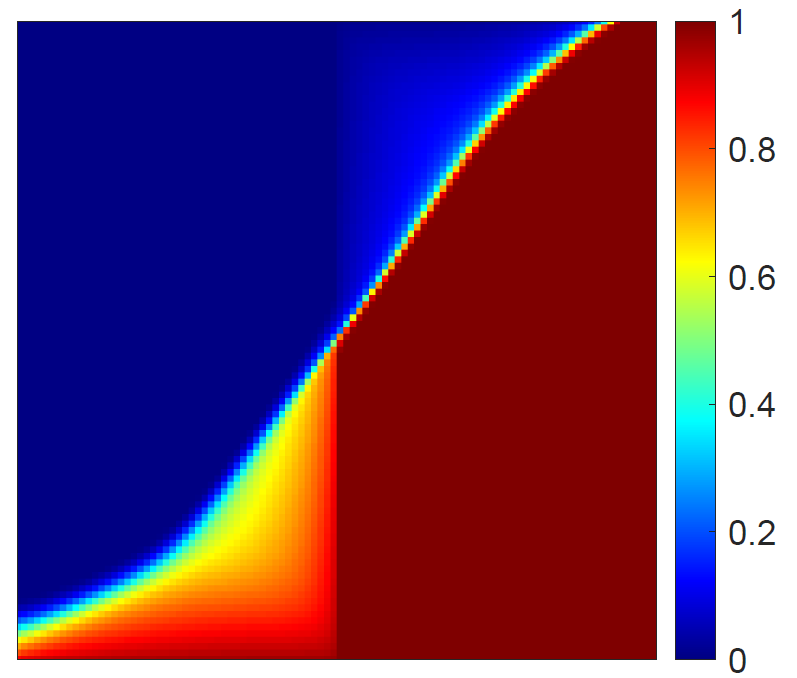}
\label{fig:lockExchangeHomoRelPerm_wetting_phase_t300}
}
\hspace{1.0cm}
\subfigure[Water saturation after 450 days]{
\includegraphics[height=0.27\textwidth]{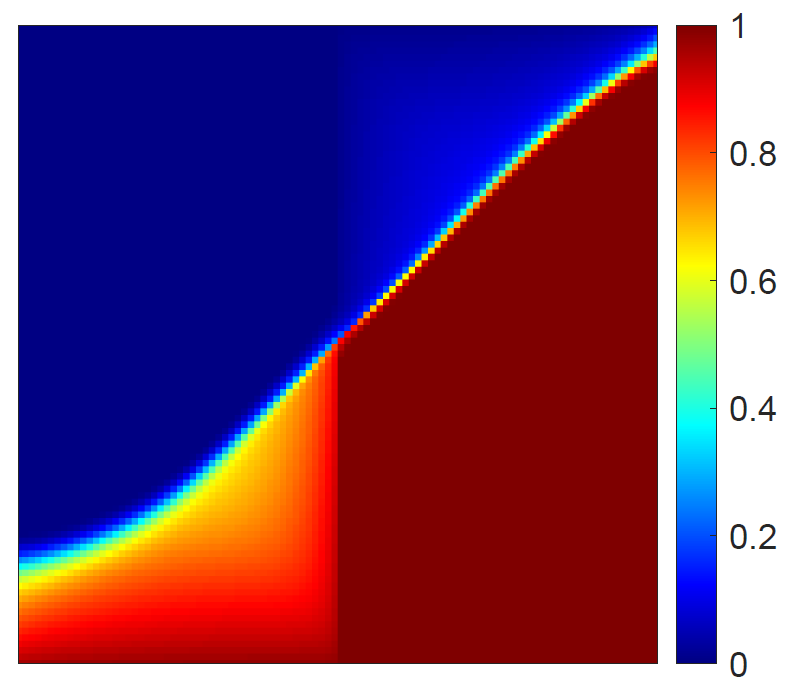}
\label{fig:lockExchangeHomoRelPerm_wetting_phase_t450}
}
\caption{Gravity segregation: water saturation maps at various times.}
\label{fig:lockExchangeHomoRelPerm_wetting_phase}
\end{figure}

We consider a two-dimensional $x$-$z$ domain of size 30.48 m $\times$ 30.48 m divided into 100 $\times$ 100 cells. 
The domain is initially saturated with the wetting phase on the left and with the non wetting phase on the right. 
The phase densities and viscosities are set to $\rho_w = 1 025$ kg.m$^{-3}$, $\rho_{nw} = 785$ kg.m$^{-3}$, $\mu_w=0.0003$ Pa.s, and $\mu_{nw} = 0.003$ Pa.s.
The phase relative permeabilities are quadratic.
The homogeneous permeability is equal to $k = 200$ mD.
We simulate 500 days of gravity segregation with a constant time step size. 
The saturation maps at different times are showed in Fig.~\ref{fig:lockExchangeHomoRelPerm_wetting_phase}.

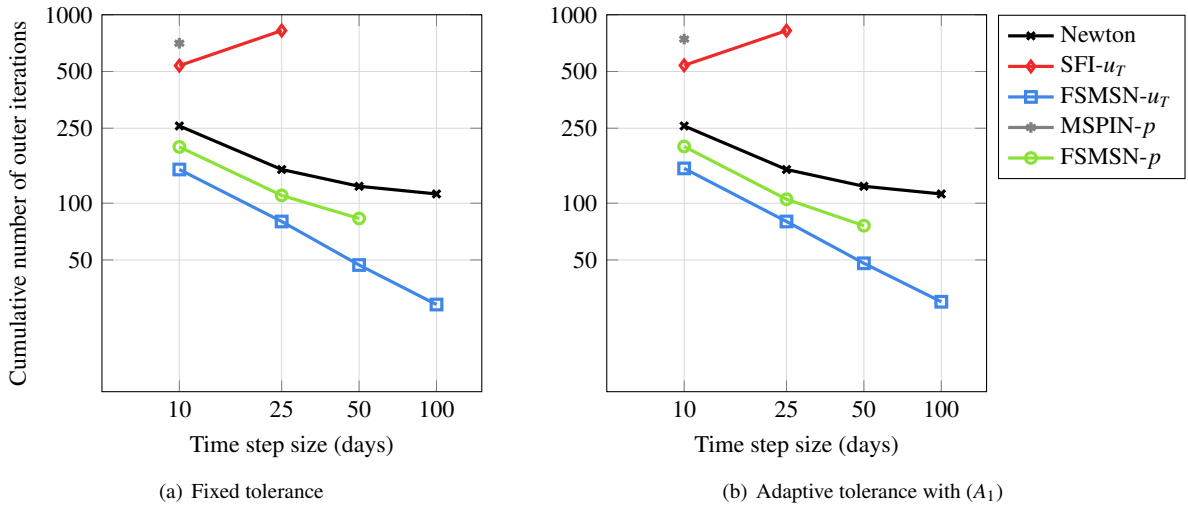
\begin{figure}[h!]
  \centering
  \subfigure[Fixed tolerance]{
  \begin{tikzpicture}
    \begin{loglogaxis}[
      width = 0.4\textwidth,
      height = 0.4\textwidth,
      grid = major,
      major grid style= {very thin,draw=gray!30},
      xmin=5,xmax=150,
      ymin=10,ymax=1000,
      xlabel={Time step size (days)},
      ylabel={Cumulative number of outer iterations},
      xtick=data,
      ytick={50,100,250,500,1000},
      yticklabels={50,100,250,500,1000},      
      ytick distance=10^1,
      ylabel near ticks,
      xlabel near ticks,
      legend style={font=\small},
      tick label style={font=\small},
      %xticklabels from table={\CFLSPETenTopLayerHighInjectionData}{cfl_labels},
      log ticks with fixed point,
      label style={font=\small},
      %legend entries={ Newton, SFI-$u_T$, FSMSN-$u_T$, MSPIN-$p$, FSMSN-$p$},
      %legend cell align={left},
      %legend pos=outer north east,
      ]
      \addplot [mark=x, black, very thick] table [x=time_dt, y=newton] {\TimestepLockExchangeData};
      \addplot [mark=diamond, scarletred1, very thick] table [x=time_dt, y=sfi_ut] {\TimestepLockExchangeData};      
      \addplot [mark=square, skyblue1, very thick] table [x=time_dt, y=fsmsn_ut] {\TimestepLockExchangeData};
      \addplot [mark=asterisk, gray, very thick] table [x=time_dt, y=mspin_p] {\TimestepLockExchangeData};      
      \addplot [mark=o, chameleon1, very thick] table [x=time_dt, y=fsmsn_p] {\TimestepLockExchangeData};
    \end{loglogaxis}
  \end{tikzpicture}
  \label{fig:nonlinearIterCurve_lockExchange_fixedTolerance}
  }\hspace{0.5cm}
  \subfigure[Adaptive tolerance with $(A_1)$]{
    \begin{tikzpicture}
    \begin{loglogaxis}[
      width = 0.4\textwidth,
      height = 0.4\textwidth,
      grid = major,
      major grid style= {very thin,draw=gray!30},
      xmin=5,xmax=150,
      ymin=10,ymax=1000,
      xlabel={Time step size (days)},
      ylabel={},
      xtick=data,
      ytick={50,100,250,500,1000},
      yticklabels={50,100,250,500,1000},      
      ylabel near ticks,
      xlabel near ticks,
      legend style={font=\small},
      tick label style={font=\small},
      log ticks with fixed point,
      label style={font=\small},
      legend entries={ Newton, SFI-$u_T$, FSMSN-$u_T$, MSPIN-$p$, FSMSN-$p$},
      legend cell align={left},
      legend pos=outer north east,
      ]
      \addplot [mark=x, black, very thick] table [x=time_dt, y=newton] {\TimestepLockExchangeAdaptivetoleranceData};
      \addplot [mark=diamond, scarletred1, very thick] table [x=time_dt, y=sfi_ut] {\TimestepLockExchangeAdaptivetoleranceData};
      \addplot [mark=square, skyblue1, very thick] table [x=time_dt, y=fsmsn_ut] {\TimestepLockExchangeAdaptivetoleranceData};
      \addplot [mark=asterisk, gray, very thick] table [x=time_dt, y=mspin_p] {\TimestepLockExchangeAdaptivetoleranceData};
      \addplot [mark=o, chameleon1, very thick] table [x=time_dt, y=fsmsn_p] {\TimestepLockExchangeAdaptivetoleranceData};  
    \end{loglogaxis}
    \end{tikzpicture}
    \label{fig:nonlinearIterCurve_lockExchange_adaptiveTolerance}
  }
  \caption{Gravity segregation: cumulative number of outer iterations for the full simulation as a function of time step size.}
  \label{fig:nonlinearIterCurve_lockExchange}
\end{figure}

The nonlinear behavior of the schemes is summarized in Fig.~\ref{fig:nonlinearIterCurve_lockExchange} and in Tables \ref{tab:lockExchange_dt10} and \ref{tab:lockExchange_dt25}.
In this challenging buoyancy-driven test case in which flow and transport are strongly coupled, SFI requires significantly more outer iterations than FIM based on Newton's method for the case with small time steps (10 days) and fails to converge when the time step is increased.
The nonlinear behavior of the solution strategies based on nonlinear preconditioning depends heavily on the formulation.
Specifically, we note that nonlinear preconditioning approaches based on fixed pressure (MSPIN-$p$ and FSMSN-$p$) do not perform as well as FSMSN-$u_T$ which used a fixed total velocity to couple the flow and transport problems.
Figure \ref{fig:nonlinearIterCurve_lockExchange} shows that both MSPIN-$p$ and FSMSN-$p$ fail to converge beyond a certain time step size (10 days for MSPIN-$p$ and 50 days for FSMSN-$u_T$). 
Importantly, FSMSN-$u_T$ exhibits a steady reduction in the number of outer iterations as a function of time step size, which is not the case for FIM with Newton's method with which the number of iterations levels off for time step sizes larger than 50 days.

Using an adaptive tolerance in the subproblems does not alter these conclusions as shown for two different time step sizes by Tables \ref{tab:adaptiveTolerance_lockExchange_dt10} and \ref{tab:adaptiveTolerance_lockExchange_dt25}.
In FSMSN-$u_T$, the three adaptive strategies considered in this work to select the subproblem tolerance cause a slight increase in the total number of outer iterations (from 151 to 154) but achieve a large reduction in the number of subproblem iterations (from 200 to 155 for flow and from 430 to 185 for transport).

%

%
%%%%%%%%%%%%%%%%%%%%%%%%%%%%%%%%%%%%%%%%%%%%%%%%%%%% Begin fig
\pgfplotstableread
{
colnames     0	   1      2      3
time_dt	     10    25	  50    100
%cfl_labels   $10^{1.5}$   $10^2$   $10^{2.5}$   $10^{3}$
fsmsn_ut     151   80   47     29
% FSMSN_p : too much transport iterations: 13190 for dt=100 and 1222 for dt=50
fsmsn_p      199   110  83     Nan
mspin_p      706   Nan  Nan    Nan    
newton       257   151  123    112  
sfi_ut       537   824  Nan    Nan
% SFI_ut timestep cut starting at dt=50
}\TimestepLockExchangeData
\pgfplotstabletranspose[colnames from=colnames]\TimestepLockExchangeData{\TimestepLockExchangeData}

\pgfplotstableread
{
colnames     0	   1      2      3
time_dt	     10    25	  50    100
fsmsn_ut     153   80     48    30 
% FSMSN_p : too much transport iterations for dt=100. for dt=50, 391 iterations
fsmsn_p      200   105    76    Nan
mspin_p      743   Nan    Nan   Nan   
newton       257   151    123   112
sfi_ut       539   824    Nan   Nan
% SFI_ut timestep cut starting at dt=50
}\TimestepLockExchangeAdaptivetoleranceData
\pgfplotstabletranspose[colnames from=colnames]\TimestepLockExchangeAdaptivetoleranceData{\TimestepLockExchangeAdaptivetoleranceData}

%%%%%%%%%%%%%%%%%%%%%%%%%%%%%%%%%%%%%%%%%%%%%%%%% End Fig
%
%
%
\begin{table}[h!]
  \centering
  \small
  \begingroup
  \setlength{\tabcolsep}{6pt}
  \begin{tabular}{l
                  c
                  c
                  c
                  c
                  c
                  c
                  c
                 }           
    \toprule
    Solver &                            && \multicolumn{2}{c}{Fixed total velocity}                        && \multicolumn{2}{c}{Fixed pressure} \\
    \cmidrule{4-5} \cmidrule{7-8} 
           & \multicolumn{1}{c}{Newton} && \multicolumn{1}{c}{SFI-$u_T$} & \multicolumn{1}{c}{FSMSN-$u_T$} && \multicolumn{1}{c}{MSPIN-$p$} & \multicolumn{1}{c}{FSMSN-$p$} \\
    \midrule
    Nonlinear iterations      & 257                   && 537   & 151  && 706   & 199 \\
    Iterations per time step  & 5.14                  && 10.74 & 3.02 && 14.12 & 3.98 \\
    Pressure iterations       & \multicolumn{1}{c}{-} && 678   & 200  && 852   & 272 \\
    Transport iterations      & \multicolumn{1}{c}{-} && 946   & 430  && 1857  & 642 \\
    \bottomrule
  \end{tabular}
  \endgroup
  \caption{Gravity segregation: nonlinear behavior of the schemes with a fixed subproblem tolerance of $10^{-6}$. We simulate 500 days with 50 time steps.}
  \label{tab:lockExchange_dt10}
\end{table}
\begin{table}[h!]
  \centering
  \small
  \begingroup
  \setlength{\tabcolsep}{6pt}
  \begin{tabular}{l      j      j      j      c      j      j      j      c      j      j      j      c      j      j      j} 
    \toprule
    Solver                    & \multicolumn{7}{c}{Fixed total velocity}                         && \multicolumn{7}{c}{Fixed pressure} \\
    \cmidrule{2-8} \cmidrule{10-16} 
                              & \multicolumn{3}{c}{SFI-$u_T$} && \multicolumn{3}{c}{FSMSN-$u_T$} && \multicolumn{3}{c}{MSPIN-$p$} && \multicolumn{3}{c}{FSMSN-$p$} \\
    \cmidrule{2-4}  \cmidrule{6-8}  \cmidrule{10-12}  \cmidrule{14-16}
                              & $A_1$ & $A_2$ & $A_3$   && $A_1$ & $A_2$ & $A_3$     && $A_1$ & $A_2$ & $A_3$   && $A_1$ & $A_2$ & $A_3$ \\
    \midrule
    Nonlinear iterations      & 539   & 558   & 558   && 153  & 154  & 154  && 743   & 717  & 712   && 200 & 200 & 200 \\
    Iterations per time step  & 10.78 & 11.16 & 11.16 && 3.06 & 3.08 & 3.08 && 14.86 & 14.34 & 14.24 && 4 & 4 & 4 \\
    Pressure iterations       & 567   & 580   & 559   && 155  & 155  & 155  && 801   & 1859 & 2909  && 203 & 201 & 204 \\
    Transport iterations      & 572   & 592   & 590   && 186  & 185  & 185  && 10761 & 10286 & 9956  && 262 & 208 & 218 \\
    \bottomrule
  \end{tabular}
  \endgroup
  \caption{Gravity segregation: nonlinear behavior of the schemes with an adaptive subproblem tolerance computed using the strategies $(A_1)$, $(A_2)$, and $(A_3)$ of Section~\ref{sec:adaptive_tolerance}. We simulate 500 days with 50 time steps.} 
\label{tab:adaptiveTolerance_lockExchange_dt10}
\end{table}
\begin{table}[h!]
  \centering
  \small
  \begingroup
  \setlength{\tabcolsep}{6pt}
  \begin{tabular}{l
                  c
                  c
                  c
                  c
                  c
                  c
                  c
                 }           
    \toprule
    Solver &                            && \multicolumn{2}{c}{Fixed total velocity}                        && \multicolumn{2}{c}{Fixed pressure} \\
    \cmidrule{4-5} \cmidrule{7-8} 
           & \multicolumn{1}{c}{Newton} && \multicolumn{1}{c}{SFI-$u_T$} & \multicolumn{1}{c}{FSMSN-$u_T$} && \multicolumn{1}{c}{MSPIN-$p$} & \multicolumn{1}{c}{FSMSN-$p$} \\
    \midrule
    Nonlinear iterations      & 151                   && 824  & 80  && d.n.c. & 110 \\
    Iterations per time step  & 7.55                  && 41.2 & 4   && d.n.c. & 5.5 \\
    Pressure iterations       & \multicolumn{1}{c}{-} && 940  & 115 && d.n.c. & 190 \\
    Transport iterations      & \multicolumn{1}{c}{-} && 1378 & 265 && d.n.c. & 509 \\
    \bottomrule
  \end{tabular}
  \endgroup
  \caption{Gravity segregation: nonlinear behavior of the schemes with a fixed subproblem tolerance of $10^{-6}$. We simulate 500 days with 20 time steps. d.n.c. denotes lack of convergence.}
  \label{tab:lockExchange_dt25}
\end{table}
\begin{table}[h!]
  \centering
  \small
  \begingroup
  \setlength{\tabcolsep}{6pt}
  \begin{tabular}{l      j      j      j      c      j      j      j      c      j      j      j      c      j      j      j} 
    \toprule
    Solver                    & \multicolumn{7}{c}{Fixed total velocity}                         && \multicolumn{7}{c}{Fixed pressure} \\
    \cmidrule{2-8} \cmidrule{10-16} 
                              & \multicolumn{3}{c}{SFI-$u_T$} && \multicolumn{3}{c}{FSMSN-$u_T$} && \multicolumn{3}{c}{MSPIN-$p$} && \multicolumn{3}{c}{FSMSN-$p$} \\
    \cmidrule{2-4}  \cmidrule{6-8}  \cmidrule{10-12}  \cmidrule{14-16}
                              & $A_1$ & $A_2$ & $A_3$   && $A_1$ & $A_2$ & $A_3$     && $A_1$ & $A_2$ & $A_3$   && $A_1$ & $A_2$ & $A_3$ \\
    \midrule
    Nonlinear iterations      & 824 & 841    & 842  && 80  & 82  & 82  && d.n.c. & d.n.c. & d.n.c. && 105 & 105 & 105 \\
    Iterations per time step  & 41.2 & 42.05 & 42.1 && 4   & 4.1 & 4.1 && d.n.c. & d.n.c. & d.n.c. && 5.25 & 5.25 & 5.25 \\
    Pressure iterations       & 844  & 8549  & 845  && 87  & 83  & 83  && d.n.c. & d.n.c. & d.n.c. && 107 & 106 & 107 \\
    Transport iterations      & 880  & 8570  & 901  && 146 & 135 & 134 && d.n.c. & d.n.c. & d.n.c. && 247 & 300 & 231 \\
    \bottomrule
  \end{tabular}
  \endgroup
  \caption{Gravity segregation: nonlinear behavior of the schemes with an adaptive subproblem tolerance computed using the strategies $(A_1)$, $(A_2)$, and $(A_3)$ of Section~\ref{sec:adaptive_tolerance}. We simulate 500 days with 20 time steps. d.n.c. denotes lack of convergence.}
\label{tab:adaptiveTolerance_lockExchange_dt25}
\end{table}

\newpage
\subsection{Fractured heterogeneous two-dimensional model}

\begin{figure}[!h]
\centering
  \subfigure[Permeability map]{
    \centering
    \includegraphics[width=0.27\textwidth]{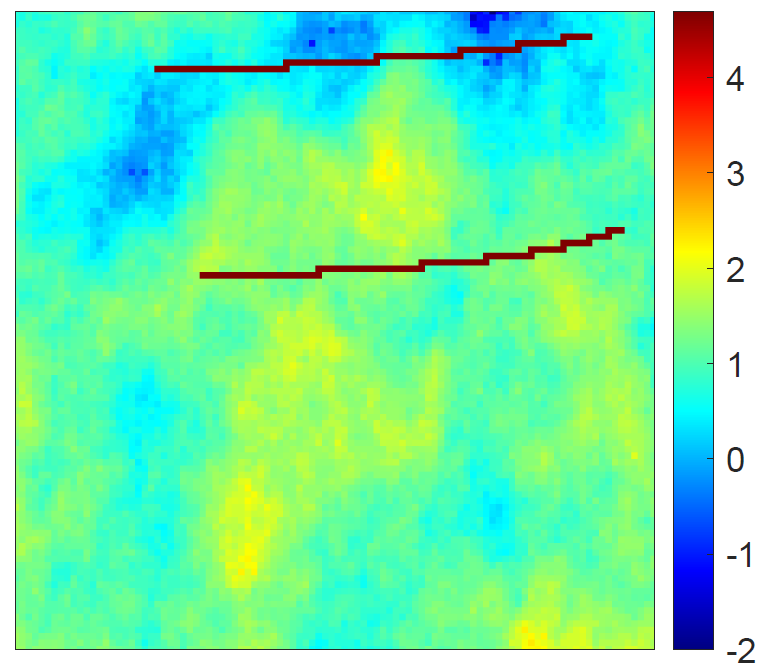}
    \label{fig:frac_perm_map}
  }
  \hfill
  \subfigure[Water saturation after 300 days.]{
    \includegraphics[width=0.27\textwidth]{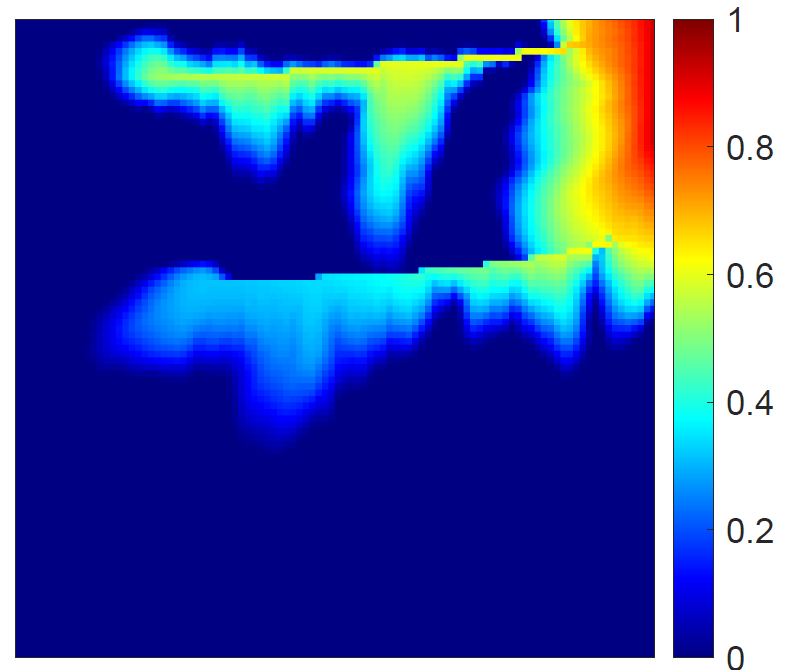}
    \label{fig:frac_sat_map_T300}
  }
  \hfill
  \subfigure[Water saturation after 600 days.]{
    \includegraphics[width=0.27\textwidth]{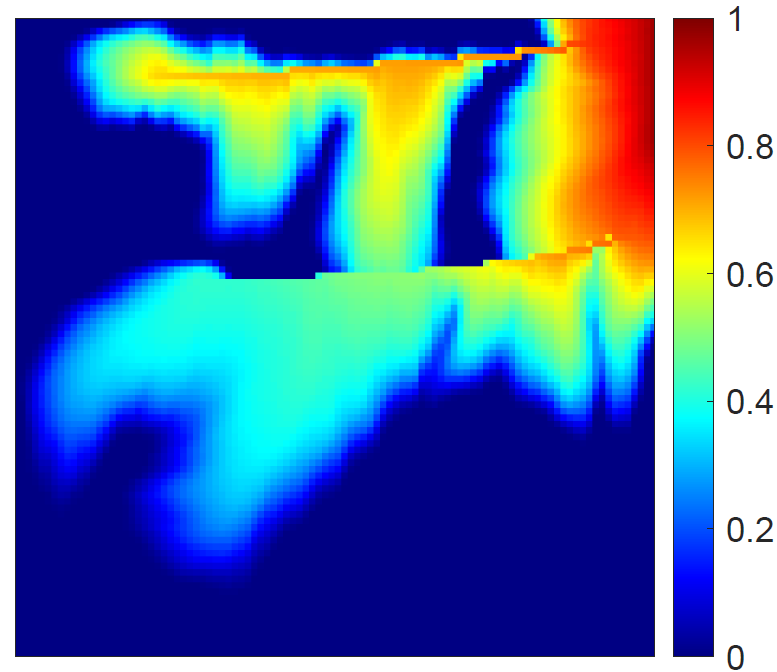}
    \label{fig:frac_sat_map_T600}
  }
  \hfill
  \subfigure[Water saturation after 900 days.]{
    \includegraphics[width=0.27\textwidth]{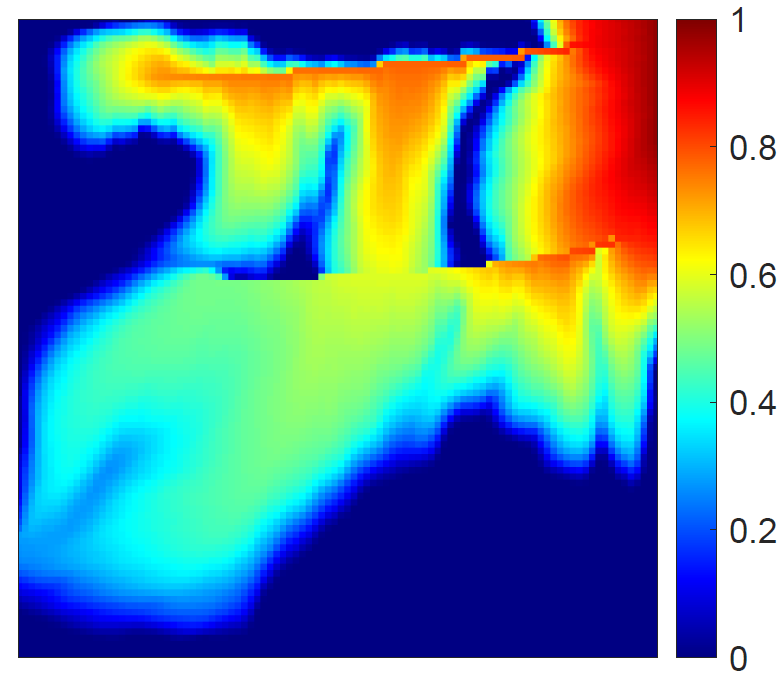}
    \label{fig:frac_sat_map_T900}
  }
  \hfill
  \subfigure[Water saturation after 1,200 days.]{
    \includegraphics[width=0.27\textwidth]{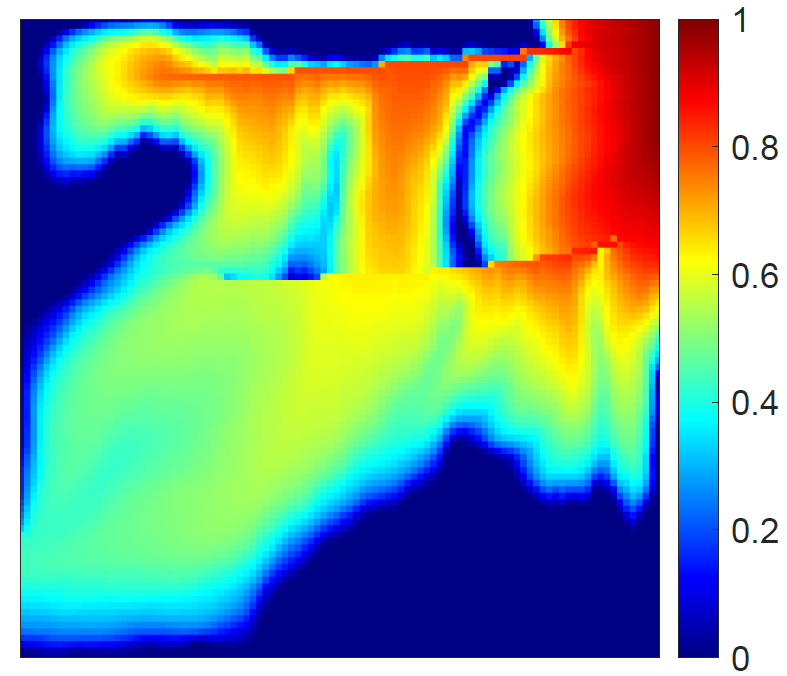}
    \label{fig:frac_sat_map_T1500}
  }
  \hfill
  \subfigure[Water saturation after 1,500 days.]{
    \includegraphics[width=0.27\textwidth]{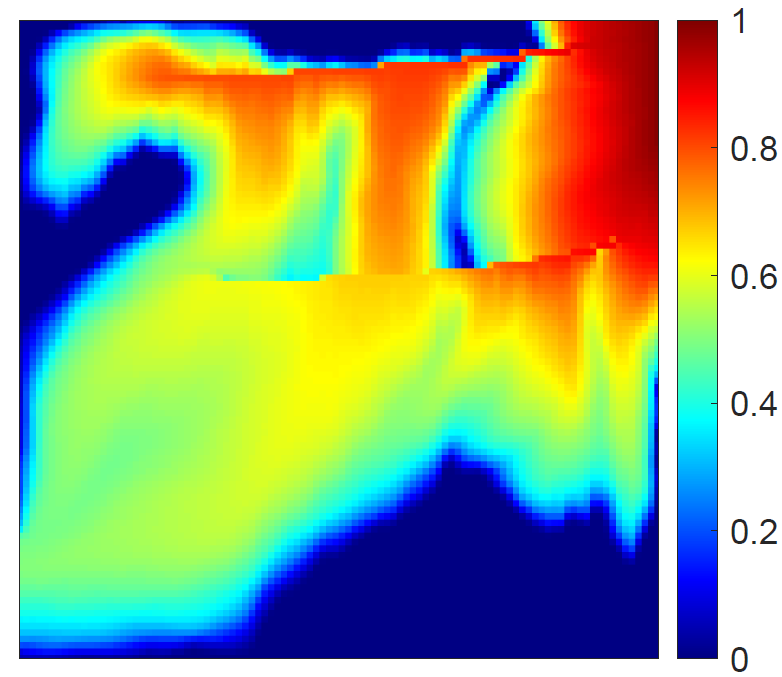}
    \label{fig:frac_sat_map_T1200}
  }
  \caption{Permeability map of the fractured heterogeneous 2D porous media and water saturation maps at various times.}
  \label{fig:fracturedMedia2D}
\end{figure}

To conclude this study, we construct a test case in which the flow is driven by competing viscous and buoyancy forces.
We consider a two-dimensional $x$-$z$ domain with a channelized permeability field--the channels can be seen as fractures modeled by contiguous cells with a very large permeability.
The fluid properties are the same as in the previous example. 
The domain is initially saturated with the non-wetting phase.
The wetting phase is injected through a well perforating twenty cells in the top-right part of the domain, while a producer perforates twenty cells in the bottom-left part of the domain. 
We simulate 1,500 days of injection (0.56 total pore volume injected). 
The permeability maps as well as the saturation map at different times is shown in Fig.~\ref{fig:fracturedMedia2D}

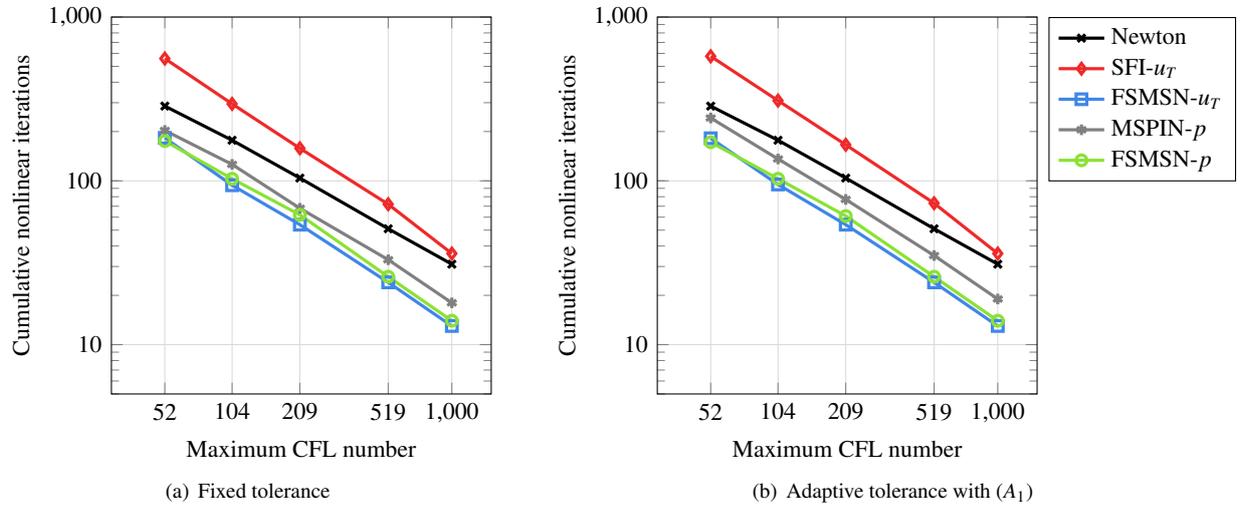
\begin{figure}[h!]
  \centering
  \subfigure[Fixed tolerance]{
  \begin{tikzpicture}
    \begin{loglogaxis}[
      width = 0.4\textwidth,
      height = 0.4\textwidth,
      grid = major,
      major grid style= {very thin,draw=gray!30},
      xmin=30,xmax=1500,
      ymin=5,ymax=1000,
      xlabel={Maximum CFL number},
      ylabel={Cumulative nonlinear iterations},
      xtick=data,
      ytick distance=10^1,
      ylabel near ticks,
      xlabel near ticks,
      legend style={font=\small},
      tick label style={font=\small},
      log ticks with fixed point,
      label style={font=\small},
      ]
      \addplot [mark=x, black, very thick] table [x=cfl, y=newton] {\CFLFracHeteroData};
      \addplot [mark=diamond, scarletred1, very thick] table [x=cfl, y=sfi_ut] {\CFLFracHeteroData};
      \addplot [mark=square, skyblue1, very thick] table [x=cfl, y=fsmsn_ut] {\CFLFracHeteroData};
      \addplot [mark=asterisk, gray, very thick] table [x=cfl, y=mspin_p] {\CFLFracHeteroData};      
      \addplot [mark=o, chameleon1, very thick] table [x=cfl, y=fsmsn_p] {\CFLFracHeteroData};
    \end{loglogaxis}
  \end{tikzpicture}
  }\hspace{0.5cm}
  \subfigure[Adaptive tolerance with $(A_1)$]{
  \begin{tikzpicture}
    \begin{loglogaxis}[
      width = 0.4\textwidth,
      height = 0.4\textwidth,
      grid = major,
      major grid style= {very thin,draw=gray!30},
      xmin=30,xmax=1500,
      ymin=5,ymax=1000,
      xlabel={Maximum CFL number},
      ylabel={Cumulative nonlinear iterations},
      xtick=data,
      ytick distance=10^1,
      ylabel near ticks,
      xlabel near ticks,
      legend style={font=\small},
      tick label style={font=\small},
      log ticks with fixed point,
      label style={font=\small},
      legend entries={ Newton, SFI-$u_T$, FSMSN-$u_T$, MSPIN-$p$, FSMSN-$p$},
      legend cell align={left},
      legend pos=outer north east,
      ]
      \addplot [mark=x, black, very thick] table [x=cfl, y=newton] {\CFLFracHeteroAdaptiveToleranceData};
      \addplot [mark=diamond, scarletred1, very thick] table [x=cfl, y=sfi_ut] {\CFLFracHeteroAdaptiveToleranceData};
      \addplot [mark=square, skyblue1, very thick] table [x=cfl, y=fsmsn_ut] {\CFLFracHeteroAdaptiveToleranceData};
      \addplot [mark=asterisk, gray, very thick] table [x=cfl, y=mspin_p] {\CFLFracHeteroAdaptiveToleranceData};      
      \addplot [mark=o, chameleon1, very thick] table [x=cfl, y=fsmsn_p] {\CFLFracHeteroAdaptiveToleranceData};
    \end{loglogaxis}
  \end{tikzpicture}
  }
  \caption{Fractured heterogeneous model: cumulative number of the nonlinear iterations as a function of the maximum CFL number observed during the simulation.}
  \label{fig:frac_nonlinear_iter_curves}
\end{figure}

The results with fixed subproblem tolerance and adaptive tolerance are in Tables \ref{tab:subtol1_frac_dt100} and \ref{tab:adaptiveTolerance_2D_fractured_media_cfl209}, respectively.
The sensitivity of the nonlinear behavior of the schemes to the time step (measured by CFL number) is shown in Fig.~\ref{fig:frac_nonlinear_iter_curves}.
Although all the solution strategies converge well in this case, the results confirm the observations made in Sections~\ref{subsec:spe10_bottom_layer} and \ref{subsec:spe10_top_layer}.
In particular, comparing the slopes of the curves for Newton-based FIM and the nonlinear preconditioners in Fig.~\ref{fig:frac_nonlinear_iter_curves} is very insightful, as it shows the improved robustness of FSMSN and MSPIN for large time step sizes. 

%
%
%Lx = 1000 ft, Ly = 1 ft, Lz = 1000 ftNx = 100, Ny = 1, Nz = 100++++++++++++ Rock properties ++++++++++++ \n
%Heterogeneous perm -> Perm filename:heterogeneous_case_100_100_1.txt
% Minimum pore volume for active blocks: PV_min = 0.000000000100000 bbl\n ++++++++++++ Fluid properties ++++++++++++ \n
%Initialization:OilSaturated
% Initial pressure at top: p = 6000 psi
%FluidProps filename : \nfluidProps.txt
%Density: rho_o = 49.0 lbm.ft^-3, rho_w = 64.0 lbm.ft^-3
%RockFluidProps filename : \nrockFluidProps.txt
%Swc = 0.00, Sor = 0.00, Sorm = 0.50Capillary pressure: No
%No hysteresis on relative permeabilities
%No hysteresis on capillary pressureS
%\n ++++++++++++ Wells ++++++++++++ \n
%Injection wells: 20 Production wells: 20
%Injectors rate constraint: q_i = 1.00 bbl.day^-1 Injectors BHP constraint: BHP = 100000.00 psi Injectors control: 0
%Days between injections: Inf, days of injection: Inf, shut all wells at time: InfProducers rate constraint: q_p = 5.00 bbl.day^-1 Producers BHP constraint: BHP = 4000.00 psi Producers control\
%: 2
%max CFL = 54.471773

%

%In Fig~\ref{fig:frac_nonlinear_iter_curves}, we show the sensitivity of the nonlinear behavior to the timestep size. The number of outer nonlinear iterations decrease as the timestep increase (even for large CFL).
\pgfplotstableread
{
colnames     0	   1      2      3     4
cfl	     52    104	  209    519   1000
%cfl_labels   $10^{1.5}$   $10^2$   $10^{2.5}$   $10^{3}$
fsmsn_ut     183   94     54     24    13
fsmsn_p      175   103    62     26    14
mspin_p      203   126    68     33    18
newton       286   177    104    51    31
sfi_ut       558   295    158    72    36
}\CFLFracHeteroData
\pgfplotstabletranspose[colnames from=colnames]\CFLFracHeteroData{\CFLFracHeteroData}
\pgfplotstableread
{
colnames     0	   1      2      3     4
cfl	     52    104	  209    519   1000
%cfl_labels   $10^{1.5}$   $10^2$   $10^{2.5}$   $10^{3}$
fsmsn_ut     182   95     54     24    13
fsmsn_p      172   103    61     26    14
mspin_p      242   136    77     35    19
newton       286   177    104    51    31
sfi_ut       576   309    166    73    36
}\CFLFracHeteroAdaptiveToleranceData
\pgfplotstabletranspose[colnames from=colnames]\CFLFracHeteroAdaptiveToleranceData{\CFLFracHeteroAdaptiveToleranceData}

%We present in Table \ref{tab:subtol1_frac_dt100} the comparison results of the different algorithms. 
\begin{table}[h!]
  \centering
  \small
  \begingroup
  \setlength{\tabcolsep}{6pt}
  \begin{tabular}{l     c      c      c     c     c     c     c}           
    \toprule
    Solver &                            && \multicolumn{2}{c}{Fixed total velocity}                        && \multicolumn{2}{c}{Fixed pressure} \\
    \cmidrule{4-5} \cmidrule{7-8} 
           & \multicolumn{1}{c}{Newton} && \multicolumn{1}{c}{SFI-$u_T$} & \multicolumn{1}{c}{FSMSN-$u_T$} && \multicolumn{1}{c}{MSPIN-$p$} & \multicolumn{1}{c}{FSMSN-$p$} \\
    \midrule
    Nonlinear iterations      & 104                   && 158   & 54  && 68   & 62 \\
    Iterations per time step  & 6.93                  && 10.53 & 3.6 && 4.53 & 4.13 \\
    Pressure iterations       & \multicolumn{1}{c}{-} && 330   & 124 && 174  & 137 \\
    Transport iterations      & \multicolumn{1}{c}{-} && 338   & 176 && 198  & 166 \\
    \bottomrule
  \end{tabular}
  \endgroup
  \caption{Fractured heterogeneous model: nonlinear behavior of the schemes with fixed subproblem tolerance of $10^{-6}$. Approximately 2\% of the interfaces experience counter-current flow. The maximum CFL number is 209 and the total PVI is 0.56.
  }
  \label{tab:subtol1_frac_dt100}
\end{table}

%In table~\ref{tab:adaptiveTolerance_2D_fractured_media_cfl209} we show the results for the different algorithms using the adaptive tolerance approach.
\begin{table}[h!]
  \centering
  \small
  \begingroup
  \setlength{\tabcolsep}{6pt}
  \begin{tabular}{l      j      j      j      c      j      j      j      c      j      j      j      c      j      j      j} 
    \toprule
    Solver                    & \multicolumn{7}{c}{Fixed total velocity}                         && \multicolumn{7}{c}{Fixed pressure} \\
    \cmidrule{2-8} \cmidrule{10-16} 
                              & \multicolumn{3}{c}{SFI-$u_T$} && \multicolumn{3}{c}{FSMSN-$u_T$} && \multicolumn{3}{c}{MSPIN-$p$} && \multicolumn{3}{c}{FSMSN-$p$} \\
    \cmidrule{2-4}  \cmidrule{6-8}  \cmidrule{10-12}  \cmidrule{14-16}
                              & $A_1$ & $A_2$ & $A_3$   && $A_1$ & $A_2$ & $A_3$     && $A_1$ & $A_2$ & $A_3$   && $A_1$ & $A_2$ & $A_3$ \\
    \midrule
    Nonlinear iterations      & 166 & 166  & 166 && 54 & 54 & 54 && 77  & 83  & 78  && 61 & 61 & 61 \\
    Pressure iterations       & 168 & 1356 & 168 && 56 & 56 & 56 && 82  & 87  & 83  && 64 & 64 & 64 \\
    Transport iterations      & 203 & 203  & 203 && 92 & 91 & 92 && 126 & 123 & 123 && 96 & 96 & 96 \\
    \bottomrule
  \end{tabular}
  \endgroup
  \caption{Fractured heterogeneous model: nonlinear behavior of the schemes with an adaptive tolerance computed using the strategies $(A_1)$, $(A_2)$, and $(A_3)$ of Section~\ref{sec:adaptive_tolerance}. The maximum CFL number is 209 and the total PVI is 0.56.} 
\label{tab:adaptiveTolerance_2D_fractured_media_cfl209}
\end{table}

%%%%%%%%%%%%%%%%%%%%%%%%%%%%%%%%%%%%%%%%%%%%%%%%%%%%%%%%%%%%%%%%%%%%%%%%%%%%%%%%%%%%
%                           Conclusion
%%%%%%%%%%%%%%%%%%%%%%%%%%%%%%%%%%%%%%%%%%%%%%%%%%%%%%%%%%%%%%%%%%%%%%%%%%%%%%%%%%%%

\section{Conclusion}
\label{sec:conclude}

Solving the nonlinear systems that result from a fully implicit discretization of the PDEs governing multiphase flow and transport in porous media is challenging.
To address this issue, we propose a field-split preconditioner referred to as Field-Split Multiplicative Schwarz Newton (FSMSN). 
The FSMSN-preconditioned iteration relies on two steps: a preconditioning step in which we solve sequentially a flow problem followed by a transport problem using a loose nonlinear tolerance; a global step in which we compute a Newton update for pressure and saturations by linearizing the preconditioned system. 
We compare its nonlinear behavior to another preconditioner of the same class (Multiplicative Schwarz Preconditioned Inexact Newton) and to standard solution strategies like Newton's method with damping to the full system, and the Sequential Fully Implicit method.

The numerical examples show that FSMSN can successfully reduce the number of outer iterations for challenging viscous-dominated and gravity-dominated problems, compared to the other solution strategies considered here.
Our results also demonstrate that this robust nonlinear behavior is preserved for large CFL numbers (corresponding to large time steps).
This is key to make sure that time step sizes can be chosen based on accuracy considerations and not constrained by the nonlinear behavior of the solution strategy.

Two key steps have to be taken to show that the improved nonlinear behavior of FSMSN can result in a reduction of the computational cost of the simulation (i.e., reduction in wall-clock time).
First, we plan to design a more adaptive FSMSN in which the preconditioning step would be used only when necessary, that is, for the first outer iterations of time steps with bad initial guesses and/or large sizes, when the large reduction in the number of outer iterations obtained with FSMSN is more likely to offset the overhead caused by the preconditioning
step.
As soon as the state of the system is close to the solution, the global step would be sufficient to enter the quadratic convergence regime.
Second, we will substitute the direct solvers used in this work with inexact Krylov-type iterative solution strategies \cite{GMRES_Saad_1986,HybridKrylovMethods_BrwonSaad_1990,AcceleratedInexatNewtonSchemesForLargeSystemsOfNonlinearEquations_FokkemaSleijpenVanderVorst}. This will reduce the overhead caused by the preconditioning step by leveraging the efficiency of specialized solvers for the subproblems, like Algebraic MultiGrid (AMG) for the pressure problem and an optimized ordering-based solver for transport \cite{PotentialBasedReducedNewtonAlgorithmForNonlinearMultiphaseFlowInPorousMedia_KwokTchelepi}. Switching to iterative solvers will also enable the use of an adaptive linear tolerance \cite{eisenstat94} in these subproblems. These improvements will enhance the efficiency of FSMSN without compromising its robust nonlinear behavior.

\section*{Acknowledgments}
Funding was provided by TotalEnergies through the FC-MAELSTROM project. The authors thank the SUPRI-B affiliates program at Stanford University and well as Joshua A. White, Nicola Castelletto (Lawrence Livermore National Laboratory), and Herv\'e Gross (TotalEnergies) for their insight and guidance.

\section*{References}
\bibliography{references}
%%%%%%%%%%%%%%%%%%%%%%%%%%%%%%%%%%%%%%%%%%%%%%%%%%%%%%%%%%%%%%% Appendix
\appendix

%%%%%%%%%%%%%%%%%%%%%%%%%%%%%%%%%%%%%%%%%%%%%%%%%%%%%%%%%%%%%%% 

\end{document}